\documentclass[mathpazo]{cicp}

\usepackage{layouts}
\usepackage{printlen}
\usepackage{algorithm}
\usepackage[noend]{algpseudocode}
\usepackage{dirtytalk}
\usepackage{lipsum}
\usepackage{amsfonts}
\usepackage{graphicx}
\usepackage{epstopdf}
\usepackage{epsfig}
\usepackage{amssymb,amsmath,amsthm}

\usepackage{amsmath,bm}
\usepackage{tikz}
\usepackage{comment}
\usepackage[mathscr]{euscript}
\usepackage{algorithm}
\usepackage{caption}
\usepackage{csquotes}
\usepackage{graphics}
\usepackage{tabularx}
\usepackage{subcaption}
\usepackage{multirow}
\usepackage{graphicx}
\usepackage{float}
\usetikzlibrary{patterns}
\usepackage{pgfplots}

\newcommand{\bkt}[1]{\left(#1\right)}

\usepackage{hyperref}

\DeclareSymbolFont{rsfs}{U}{rsfs}{m}{n}
\DeclareSymbolFontAlphabet{\mathscrsfs}{rsfs}
\numberwithin{equation}{section}

\ifpdf
  \DeclareGraphicsExtensions{.eps,.pdf,.png,.jpg}
\else
  \DeclareGraphicsExtensions{.eps}
\fi

\makeatletter
\def\imod#1{\allowbreak\mkern10mu({\operator@font mod}\,\,#1)}
\makeatother

\begin{document}
\title{A new Directional Algebraic Fast Multipole Method based iterative solver for the Lippmann-Schwinger equation accelerated with HODLR preconditioner}


\author[Vaishnavi G et.~al.]{Vaishnavi Gujjula\affil{1}\comma\corrauth
       and Sivaram Ambikasaran\affil{1}}
\address{\affilnum{1}\ Department of Mathematics,
        Indian Institute of Technology Madras,
         Chennai 600036, India.
  }
\emails{{\tt vaishnavihp@gmail.com} (V.~Gujjula),
         {\tt sivaambi@smail.iitm.ac.in} (S.~Ambikasaran)}

\begin{abstract}
  We present a fast iterative solver for scattering problems in 2D, where a penetrable object with compact support is considered. By representing the scattered field as a volume potential in terms of the Green's function, we arrive at the Lippmann-Schwinger equation in integral form, which is then discretized using an appropriate quadrature technique.
  The discretized linear system is then solved using an iterative solver accelerated by Directional Algebraic Fast Multipole Method (DAFMM).
  The DAFMM presented here relies on the directional admissibility condition of the 2D Helmholtz kernel~\cite{engquist2009fast}, and the construction of low-rank factorizations of the appropriate low-rank matrix sub-blocks is based on our new Nested Cross Approximation (NCA)~\cite{https://doi.org/10.48550/arxiv.2203.14832}.
  The advantage of the NCA described in~\cite{https://doi.org/10.48550/arxiv.2203.14832} is that the search space of so-called far-field pivots is smaller than that of the existing NCAs~\cite{bebendorf2012constructing,zhao2019fast}.
  Another significant contribution of this work is the use of HODLR based direct solver~\cite{ambikasaran2013mathcal} as a preconditioner to further accelerate the iterative solver.
  In one of our numerical experiments, the iterative solver does not converge without a preconditioner. We show that the HODLR preconditioner is capable of solving problems that the iterative solver can not.
  Another noteworthy contribution of this article is that we perform a comparative study of the HODLR based fast direct solver, DAFMM based fast iterative solver, and HODLR preconditioned DAFMM based fast iterative solver for the discretized Lippmann-Schwinger problem. To the best of our knowledge, this work is one of the first to provide a systematic study and comparison of these different solvers for various problem sizes and contrast functions. In the spirit of reproducible computational science, the implementation of the algorithms developed in this article is made available at \url{https://github.com/vaishna77/Lippmann_Schwinger_Solver}.
\end{abstract}


\maketitle

\section{Introduction}
\label{section:Introduction}
This article focuses on developing a fast iterative solver for scattering problems in 2D.
Consider a penetrable object with an electric susceptibility (or \textit{contrast function}) of $q(x)$. Assume $q(x)$ to have compact support in a domain $\Omega$.
Let $u^{inc}(x)$ be the incident field and $u^{scat}(x)$ be the unknown scattered field. Let $\kappa$ be the wavenumber of the incident field. The total field $u$, which is the sum of incident and scattered fields, follows the time-harmonic Helmholtz equation
  \begin{equation} \label{eq:1}
    \nabla^{2}u(x) +\kappa^{2}(1+q(x))u(x) = 0, \hspace{3mm} x\in \mathbb{R}^{2}.
  \end{equation}
  The incident field satisfies the homogeneous Helmholtz equation
  \begin{equation} \label{eq:2}
    \nabla^{2}u^{inc}(x) +\kappa^{2}u^{inc}(x) = 0, \hspace{3mm} x\in \Omega.
  \end{equation}
  It follows from Eq.~\eqref{eq:1} and Eq.~\eqref{eq:2} that $u^{scat}(x)$ satisfies
  \begin{equation} \label{eq:3}
    \nabla^{2}u^{scat}(x) +\kappa^{2}(1+q(x))u^{scat}(x) = -\kappa^{2}q(x)u^{inc}(x), \hspace{3mm} x\in \Omega.
  \end{equation}

  To ensure the scattered field propagates to infinity without any spurious resonances, we enforce the Sommerfeld radiation condition
  \begin{equation} \label{eq:4}
    \lim_{r\to \infty} r^{1/2}\bkt{\frac{\partial u}{\partial r}^{scat} - i\kappa u^{scat}}=0, \hspace{3mm}\text{ where }r=||x||.
  \end{equation}
  There exist many techniques to solve the scattered field. A few of them worth mentioning are: constructing a variational form, discretizing the differential operator, reformulating it as a volume integral equation.
 We use the volume integral equation technique as described in~\cite{ambikasaran2016fast}, where the scattered field is expressed as a volume potential
  \begin{equation} \label{eq:5}
    u^{scat}(x) = V[\psi](x) = \int_{\Omega}G_{\kappa}(x,y)\psi(y)dy, \hspace{3mm} x\in \Omega
 \end{equation}
 where
 \begin{equation}
   G_{\kappa}(x,y) = \frac{i}{4}H_{0}^{(1)}(\kappa ||x-y||)
\end{equation}
is the Green's function of Helmholtz equation in 2D.
 Using Eq.~\eqref{eq:3} and Eq.~\eqref{eq:5}, we obtain the Lippmann-Schwinger equation
 \begin{equation} \label{eq:6}
   \psi(x) + \kappa^{2}q(x)V[\psi](x) = f(x), \hspace{3mm} x\in \Omega
 \end{equation}
 where $f(x)=-\kappa^{2}q(x)u^{inc}(x)$.
 The task is to numerically solve for $\psi(x)$ and then obtain $u^{scat}(x)$.

The present article discusses a fast iterative solver for the Lippmann-Schwinger equation, developed on an adaptive grid. Iterative solvers rely on matrix-vector products, which can be prohibitively expensive for large-sized problems when the underlying matrix is dense. It scales as $\mathcal{O}(N^2)$, where $N$ is the number of unknowns in the discretized form of the above linear system. Our solver exploits the low-rank sub-blocks in the underlying matrix, which reduces the complexity of the algorithm to $\mathcal{O}(N\log(N))$.

Matrix-vector product can be interpreted as an $N$-body problem. Such $N$-body problems arising from the Helmholtz kernel have been studied extensively in the literature.

Rokhlin (1990)~\cite{Rokhlin1990} developed the high frequency fast multipole method (HF-FMM) wherein the radiation fields are expressed as partial wave expansions and diagonal translation operators were developed.

Engquist and Ying (2009)~\cite{engquist2009fast} developed a directional algorithm for an $N$-body problem. It relies on the \textit{directional} admissibility condition to identify low-rank matrix sub-blocks.
With the directional admissibility condition, the rank of the low-rank sub-blocks is independent of wavenumber.
The matrix sub-blocks that can be low rank approximated were compressed using the pivoted QR factorization. And the
pivots were chosen using a randomized sampling technique.

Messner et al. (2012)~\cite{messner2012fast} too used the directional admissibility condition for low-rank and developed a Chebyshev interpolation-based summation technique.

  Directional $\mathcal{H}^{2}$ matrices were introduced by Bebendorf et al. (2015) in~\cite{bebendorf2015wideband}, which is a sub-class of $\mathcal{H}^{2}$ matrices~\cite{borm2003introduction,borm2003hierarchical,hackbusch2015mathcal,hackbusch2002h2,inbook} with the directional admissibility condition of Helmholtz kernel.
  Using the Directional $\mathcal{H}^{2}$ structure in the high frequency regime, an algebraic summation technique for a 3D Helmholtz integral operator, arising out of Galerkin discretization was presented.
  The bases vectors of the matrix sub-blocks that can be low rank approximated were formulated using nested cross approximation (NCA)~\cite{bebendorf2012constructing} -
  a technique that develops nested bases.
  B{\"o}rm (2017)~\cite{borm2017directional} used Directional $\mathcal{H}^{2}$ matrices and developed a summation technique for a 3D Helmholtz integral operator, arising out of Galerkin discretization. The low-rank compressions were formulated using an adaptive QR factorization.

  In this article we will be looking at a fast, directional, and algebraic method.
  The low-rank approximations in this article are constructed based on our new NCA~\cite{https://doi.org/10.48550/arxiv.2203.14832}, developed for a sub-class of $\mathcal{H}^{2}$ matrices which we term as \textit{FMM matrices}, that follow strong admissibility condition, i.e., the interaction between neighboring cluster of particles is considered full-rank and the interaction between \textit{well-separated} cluster of particles is approximated to low-rank.

  NCA was first introduced in~\cite{bebendorf2012constructing} for non-oscillatory kernels to construct $\mathcal{H}^{2}$ matrix representation. It relies on a geometrical method to find pivots in $\mathcal{O}(N\log N)$ complexity. It was extended to the Helmholtz kernel in~\cite{bebendorf2015wideband} by constructing \textit{directional low-rank approximations}.
  An $\mathcal{O}(N)$ NCA, that is completely algebraic is developed in~\cite{zhao2019fast}.
  The NCA developed by us in~\cite{https://doi.org/10.48550/arxiv.2203.14832} differs from the existing NCAs~\cite{bebendorf2012constructing,zhao2019fast} in the technique of choosing pivots, a key step in the approximation.
  In articles~\cite{bebendorf2012constructing,zhao2019fast} the search space of far-field pivots of a cluster of points is considered to be the entire far-field region of the domain containing the support of the cluster of points. Whereas in our NCA~\cite{https://doi.org/10.48550/arxiv.2203.14832}, the search space of far-field pivots of a box of the FMM tree is limited to the union of boxes in its interaction list, an efficient representation of its far-field for FMM matrices.
  As a consequence, the time taken to construct the FMM matrix representation using our method~\cite{https://doi.org/10.48550/arxiv.2203.14832} is lesser than that of~\cite{bebendorf2012constructing,zhao2019fast}.
  We refer the readers to~\cite{https://doi.org/10.48550/arxiv.2203.14832} for numerical evidence of the accuracy of the method. We discuss more on far-field pivots in Section~\ref{sec:constructionNCA}.

  The advantages of our NCA in~\cite{https://doi.org/10.48550/arxiv.2203.14832} over the NCA of~\cite{bebendorf2012constructing} are:
  \begin{itemize}
    \item
    The time complexity of finding pivots of the former is $\mathcal{O}(N)$, whereas that of the latter is $\mathcal{O}(N\log(N))$.
    \item
    The P2M/M2M and L2L/L2P translation operators can be obtained from
    the pivot-choosing routine of the former, whereas, for the latter, they need to be computed separately.
  \end{itemize}

  In this article, we adapt our NCA~\cite{https://doi.org/10.48550/arxiv.2203.14832}, developed for non-oscillatory kernels to the 2D Helmholtz kernel by constructing directional low-rank approximations. In addition to the difference in the method of choosing pivots,
  our method differs from~\cite{bebendorf2015wideband} in the construction of directional low-rank approximation. In~\cite{bebendorf2015wideband}, only the far-field pivots are directional and the self pivots are non-directional. Whereas in this article all the pivots are directional.

  Furthermore, we develop a Directional Algebraic Fast Multipole Method (from now on abbreviated as DAFMM), to compute fast matrix-vector products, that uses NCA to low-rank approximate the appropriate matrix sub-blocks. DAFMM attempts to address some of the drawbacks of the existing techniques.
  \begin{enumerate}
    \item
    We rely on the directional admissibility condition to identify low-rank sub-blocks and an algebraic method to obtain low-rank approximations. This avoids numerical instabilities associated with the wave expansions of the Helmholtz kernel.
    \item
    We use an algebraic method to find the bases for low-rank sub-blocks.
    The advantages of an algebraic method are: (i) It is highly adaptive to the problem at hand
(ii) Can be used in a completely black-box fashion
(iii) Ranks are typically lower than the analytic methods since it is highly problem-specific.
  \end{enumerate}

	We illustrate the applicability of the DAFMM algorithm by constructing a fast iterative solver for the Lippmann-Schwinger equation.

  We summarise the key aspects of this article here:
    \begin{enumerate}
      \item
      A completely algebraic fast iterative solver for the discretized Lippmann-Schwinger equation is constructed. The matrix-vector products that are encountered are accelerated by DAFMM, which relies on our new NCA.
      \item
      We use the HODLR scheme~\cite{ambikasaran2013mathcal}, a direct solver developed for matrices whose off-diagonal blocks are low-rank, as a preconditioner to speed up the convergence of the iterative solver. As will be illustrated in the numerical examples, this preconditioner is crucial to converge to the solution.
      \item
	  A comparative study of a class of fast iterative and direct solvers (particularly for the Lippmann-Schwinger equation) is presented.
    \end{enumerate}

The rest of the article is organized as follows. In Section~\ref{sec:Discretization}, we discretize the Lippmann-Schwinger equation to obtain a discrete linear system. In Section~\ref{sec:admissibilityCriteria}, we describe the admissibility condition for low-rank.
We present NCA in Section~\ref{sec:NCA} and in Section~\ref{section:DFMM} we give a detailed description of the DAFMM algorithm, which is the primal part of the iterative solver.
In Section~\ref{section:numericalResults},
we present numerical experiments and conclude with a comparative analysis of iterative and direct solvers.

\section{Discretization}
\label{sec:Discretization}
The integral equation, Eq.~\eqref{eq:6}, is discretized as in~\cite{ambikasaran2016fast} to obtain a discrete linear system.
Instead of repeating what is done in~\cite{ambikasaran2016fast}, we summarize the important details. We refer the readers to Section 2 of ~\cite{ambikasaran2016fast} for more details.

Let $\Omega$ denote a compact square domain amenable for discretization containing the support of $q(x)$.
A quad-tree is built on the domain, and a $p\times p$ tensor product Chebyshev grid is formed in each leaf node of the tree. The unknowns are the function values $\psi(x)$ evaluated at the grid points of all the leaf nodes.

\subsection{Construction of quad-tree}
\label{sec:Constructionofquad-tree}
We hierarchically subdivide the domain based on an adaptive quad-tree data structure.
Level $0$ of the tree is the domain itself. We recursively sub-divide a box $B$ at level $\nu$ in the tree into four child boxes (belonging to level $\nu+1$) under the following conditions:

\begin{enumerate}
	\item
	the contrast function is not \textit{well-resolved} in the box $B$.
  \item
	the incident field is not \textit{well-resolved} in the box $B$.
  \item
	the box $B$ lies in the high frequency regime. Section~\ref{sec:admissibilityCriteria} has more details on it.
\end{enumerate}
The contrast function being \textit{well-resolved} means that it is approximated to a user-specified tolerance. We evaluate $q(x)$ at the grid points of the box $B$ and use it to evaluate the $p^{th}$ order coefficients of the Chebyshev approximation, say $Q^{B}$. We evaluate $Q^{B}$ at the grid points of the child boxes and if it agrees with the function values to a user-specified tolerance $\epsilon_{grid}$, we stop further refinement of the box $B$. Similar procedure is followed for $u^{inc}(x)$.

In problems with a rapidly varying contrast function or incident field, it is advantageous to use an adaptive tree.
  In this article, we use a \textit{level-restricted tree}, i.e., any two boxes that share a boundary are not more than a predetermined number of levels apart.
  This helps in ensuring that we only have a few types of neighboring (see Table~\ref{table:Notations} for the definition of a neighbor of a box) interactions to consider for each leaf box. Such interactions can be precomputed and reused across boxes of the same level thereby reducing the computational overhead. To have fewer precomputations, we impose a condition that any two boxes that share a boundary are not more than one level apart, as done in~\cite{ambikasaran2016fast}. Note that given any adaptive tree, one can construct a level-restricted tree.

\subsection{The discrete linear system}
We now construct a discrete model of the Lippmann-Schwinger equation, Eq.~\eqref{eq:6}.
In the process, we approximate $\psi(x)$ to a given precision.

To build a $\psi$ that is $p$th-order accurate, we choose polynomials in two dimensions of degree less than $p$ as basis functions, wherein the support of these basis functions is limited to the box under consideration. This means we have a total of $N_{p} = p(p+1)/2$ polynomials as basis functions in two dimensions. Let $\{b_{l}(\zeta_{1},\zeta_{2})\}_{l=1}^{N_p}$ be such polynomials, which span the space of polynomials of degree less than $p$. We use Chebyshev polynomials scaled to box $[-1,1]^{2}$. The reason for scaling the basis functions to box $[-1,1]^{2}$ is to reuse the computations involving them for other boxes as well. This reduces computational overhead.
For a leaf box $B$ of width $2\beta^{B}$, centered at $(\alpha_{1}^{B}, \alpha_{2}^{B})$, $\psi(x)$ is approximated as
\begin{equation} \label{eq:7}
  \psi(x) \approx \psi^{B}(x) = \sum_{l=1}^{N_{p}} c_{l}^{B} b_{l}\left(\frac{\zeta_{1}-\alpha_{1}^{B}}{\beta^{B}}, \frac{\zeta_{2}-\alpha_{2}^{B}}{\beta^{B}}\right)\hspace{3mm}
  \forall \hspace{1mm} x=(\zeta_{1},\zeta_{2})\in B.
\end{equation}
The co-efficients $c_l^{B}$ are so chosen to match $\psi^{B}$ to $\psi$, at the tensor product Chebyshev nodes of $B$, $\{x_{i}^B\}_{i=1}^{p^2}$, of order $p \times p$.
By evaluating Eq.~\eqref{eq:7} at each $x_{i}^{B} = (\zeta_{i,1}^B, \zeta_{i,2}^B)$, the grid points of box $B$, we form a vector $\vec{\psi}^{B}=[\psi(x_{1}^B), \psi(x_{2}^B), \ldots,\psi(x_{p^{2}}^B)]^{T}$. Vector $\vec{\psi}^{B}$, is expressed in terms of vector $\vec{c}^{B} = [{c}^{B}_{1}, {c}^{B}_{2}, \ldots,{c}^{B}_{N_{p}}]^{T}$ as
\begin{equation} \label{eq:9}
  \vec{\psi}^{B} = Q\vec{c}^{B}
\end{equation}
where $Q \in \mathbb{R}^{p^2 \times N_p}$ is the interpolation matrix, whose entries are given by
\begin{equation} \label{eq:9b}
  Q_{il} = b_{l}\left(\frac{\zeta_{i,1}^B-\alpha_{1}^{B}}{\beta^{B}}, \frac{\zeta_{i,2}^B-\alpha_{2}^{B}}{\beta^{B}}\right).
\end{equation}
By taking the pseudo-inverse of $Q$, we obtain $\vec{c}^{B}$ in terms of $\vec{\psi}^{B}$
\begin{equation} \label{eq:10}
  \vec{c}^{B} = Q^{\dagger}\vec{\psi}^{B}.
\end{equation}

We define $Q^{\dagger}(l,:)$ as in MATLAB's matrix slicing notation, i.e., $Q^{\dagger}(l,:)\in \mathbb{C}^{1 \times p^{2}}$ is the $l^{th}$ row vector of matrix $Q$. Using Eq.~\eqref{eq:10} and Eq.~\eqref{eq:7}
\begin{equation} \label{eq:11}
  \psi^{B}(x) = \sum_{l=1}^{N_{p}}  Q^{\dagger}(l,:)\vec{\psi}^{B}  b_{l}\left(\frac{\zeta_{1}-\alpha_{1}^{B}}{\beta^{B}}, \frac{\zeta_{2}-\alpha_{2}^{B}}{\beta^{B}}\right)
\end{equation}

Given a quadtree subdivision of the domain, let $\mathcal{L}$ be the set of all leaf boxes. Eq.~\eqref{eq:5}, that defines the scattered field as the volume integral, can be re-written as
\begin{equation}
  u^{scat}(x) = V[\psi](x) = \sum_{B \in \mathcal{L}}\int_{\Omega \in B}G_{\kappa}(x,y)\psi(y)dy.
\end{equation}
An approximation of $V[\psi](x)$ using Eq.~\eqref{eq:11} takes the form
\begin{equation} \label{eq:11b}
  V[\psi](x) \approx \sum_{\mathcal{L}}\int_{B}G_{\kappa}(x,y)\sum_{l=1}^{N_{p}}  Q^{\dagger}(l,:)\vec{\psi}^{B}b_{l}\left(\frac{\zeta_{1}-\alpha_{1}^{B}}{\beta^{B}}, \frac{\zeta_{2}-\alpha_{2}^{B}}{\beta^{B}}\right) dy.
\end{equation}
Using the approximate $V[\psi](x)$ in Eq.~\eqref{eq:6}, we have
\begin{equation} \label{eq:11c}
  \psi(x)+\kappa^{2}q(x) \sum_{\mathcal{L}}\int_{B}G_{\kappa}(x,y)\sum_{l=1}^{N_{p}}  Q^{\dagger}(l,:)\vec{\psi}^{B}b_{l}\left(\frac{\zeta_{1}-\alpha_{1}^{B}}{\beta^{B}}, \frac{\zeta_{2}-\alpha_{2}^{B}}{\beta^{B}}\right) dy \approx f(x).
\end{equation}
Let $\mathcal{X}=\{x_{i}\}_{i=1}^{N}$ denote the union of the grid points of all the leaf nodes of the tree, where $N=p^{2}|\mathcal{L}|$.
By enforcing Eq.~\eqref{eq:11c} at the points in $\mathcal{X}$, we obtain a discrete linear system
\begin{equation} \label{eq:12a}
A\vec{\psi} = \vec{f}
\end{equation}
where $\psi_{j} = \psi(x_{j})$, $f_{j} = f(x_{j})$ and
the $(i,j)^{th}$ matrix entry of $A$ is given by

\begin{equation} \label{eq:13}
  A_{ij} = \delta_{ij} + \kappa^{2}q(x_{i})\int_{\Omega \in B}G_{\kappa}(x_i,y) \sum_{l=1}^{N_{p}}  Q^{\dagger}_{l,j'} b_{l}\left(\frac{y_{1}-\alpha_{1}^{B}}{\beta^{B}}, \frac{y_{2}-\alpha_{2}^{B}}{\beta^{B}}\right)  dy_{1} dy_{2}
\end{equation} where
\[ j'=1+(j-1) \bmod p^{2}. \]

\subsection{Field computations}
\label{section:FieldComputations}

To obtain the matrix entries in Eq.~\eqref{eq:13}, we need to compute the integrals of the form
$$\displaystyle \int_{\Omega \in B} G_{\kappa}(x_i,y)b_{l}\left(\frac{y_{1}-\alpha_{1}^{B}}{\beta^{B}}, \frac{y_{2}-\alpha_{2}^{B}}{\beta^{B}}\right)  dy_{1} dy_{2}.$$
We use an adaptive Gauss quadrature technique to evaluate these integrals.

The integrals encountered in the far-field interactions are computed by expanding the 2D Helmholtz kernel in Eq.~\eqref{eq:13} using the Graf-Addition theorem~\cite{abramowitz1964handbook,olver2010nist} as presented in~\cite{ambikasaran2016fast}. All the relevant intermediate integrals, arising out of expanding the Helmholtz kernel, are invariant to a shift in the coordinates. So it is sufficient to evaluate and tabulate them once per level of the tree. The integrals required to compute the far-field interactions are then obtained from these tabulated integrals. This enables us to reduce the computation time for obtaining the matrix entries.

The integrals encountered in the near-field interactions are computed by evaluating the 2D Helmholtz kernel in Eq.~\eqref{eq:13} using the Boost library~\cite{Boost}. With the level-restricted tree that we considered, a box at level $\nu$ can have atmost 32 neighbors (8 neighbors belonging to level $\nu$, 12 neighbors belonging to level $\nu+1$, and 12 neighbors belonging to level $\nu-1$)~\cite{ambikasaran2016fast}. Further, the near-field integral only depends on the distance between a target point and the coordinates of a box. So it is sufficient to tabulate the near-field interactions once per level of the tree and they can be re-used to obtain the near-field interactions of all leaf boxes.

\section{Admissibility condition for low-rank}
\label{sec:admissibilityCriteria}
Given the linear system of Eq.~\eqref{eq:12a}, we propose a fast method to solve it. The method relies on making use of the fact that certain sub-blocks in the matrix can be well-approximated by low-rank matrices.
These sub-blocks are identified by the admissibility condition as discussed below.
We use different admissibility conditions for boxes in low and high frequency regimes.
Given a parameter $T$, a box of width $2w$ belonging to the quad-tree is said to be in the high frequency regime if $\bkt{\kappa w}^2 > T$, else it is said to be in the low frequency regime.
We now state the admissibility conditions for low-rank in low and high frequency regimes.

\subsection{Admissibility condition for low-rank in low frequency regime}
\label{section:admissibilityCriteriaLowRank}
The interaction between boxes $X$ and $Y$ in the low frequency regime is said to be admissible for low-rank approximation if they are separated by a distance at least equal to the width of the larger box, i.e., the centers of the boxes need to be apart by at least $4w$, where $2w$ is the width of the larger box. Such boxes $X$ and $Y$ are said to be \textit{well-separated} in the low frequency regime and the corresponding matrix sub-block that represents the interaction between boxes $X$ and $Y$ is said to be \textit{admissible}.

The \textit{far-field region of $X$}, denoted by $F(X)$ is defined as the union of all boxes that are \textit{well-separated} to $X$.

\subsection{Admissibility condition for low-rank in high frequency regime}
\label{ssec:HFR_admissibility}
If the same admissibility condition of the low frequency regime is used in the high frequency regime, the rank of the approximation grows linearly with $\kappa$.
We use the directional admissibility condition proposed by Engquist et al. in~\cite{engquist2009fast}, with which the rank of the admissible blocks is independent of $\kappa$.
We summarise it here.

\begin{figure}[H]
  \begin{center}
    \includegraphics[width=0.3\linewidth]{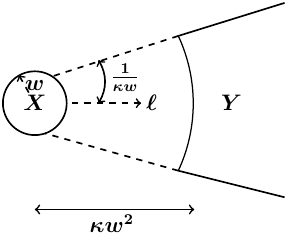}
  \end{center}
\caption{Two sets $X$ and $Y$ that satisfy the directional admissibile condition for low-rank in high frequency regime. $Y$ is separated to $X$ by a distance at least equal to ${\kappa w^{2}}$ and $Y$ lies in a cone originating from the center of $X$, which has a cone angle of ${\frac{1}{\kappa w}}$.}
\label{fig:HFR_farfield}
\end{figure}

Consider sets $X$ and $Y$, defined as
\begin{equation}\label{eq:defX}
  X = \{x:|x-x_{0}| \leq w\}
\end{equation}
\begin{equation}\label{eq:HFR_admissibility0}
  Y_{1} = \{y:|y-x_{0}| \geq \kappa w^{2}\} \hspace{5mm} and \hspace{5mm} Y_{2} = \{y:\measuredangle{(y-x_{0}, \ell)} \leq \frac{1}{\kappa w}\}
\end{equation}
\begin{equation}\label{eq:HFR_admissibility}
  Y= Y_{1}\cap Y_{2}
\end{equation}
where $\ell\in\mathbb{R}^{2}$, $x_{0}\in \mathbb{R}^{2}$ and $w>0$.
Sets $X$ and $Y$ are pictorially represented in Figure~\ref{fig:HFR_farfield}.
Then there exist functions $\{u_{i}(x)\}_{1\leq i \leq c}$ and $\{v_{i}(y)\}_{1\leq i \leq c}$ such that the 2D Helmholtz kernel $G_{\kappa}(x,y)$ has a separable expansion of the form
\begin{equation}\label{eq:separable}
	G_{\kappa}(x,y) \approx \sum\limits_{i=1}^{c} u_{i}(x)v_{i}(y),\hspace{3mm}  x\in X \text{ and } y\in Y.
\end{equation}
where $c$ is the rank of the approximation.
 For proof of Eq.~\eqref{eq:separable}, we refer the readers to~\cite{engquist2009fast}. It is to be noted that the directional admissibility condition consists of a distance and an angle condition, illustrated by $Y_{1}$ and $Y_{2}$ respectively.

Such sets $X$ and $Y$, are said to be \textit{well-separated} in the high frequency regime and the corresponding matrix sub-block representing the interaction between these two sets, is said to be \textit{admissible}, and can be well approximated by a low-rank matrix.

$Y$ is said to be the \textit{far-field region of $X$ in direction $\ell$}, denoted by $F(X,\ell)$.
\subsubsection{Construction of cones}
\label{Construction_of_cones}
In accordance with the directional admissibility condition for low-rank in the high frequency regime, the far-field region of a box is divided into multiple conical regions. We follow the procedure stated in~\cite{engquist2009fast} to construct cones, wherein
as we traverse up the tree in the high frequency regime, we divide the far-field of a box into twice the number of conical regions of its child.
We bisect a cone of a box at level $\nu+1$ to form the cones of its parent at level $\nu$.
Each conical region is indexed by its axis vector.
For a box $B$, the set of all axis vectors of the cones associated with its far-field region is denoted by $L(B)$.

\section{Nested Cross Approximation}\label{sec:NCA}
The construction of low-rank approximations of the admissible sub-blocks of the matrix is based on our new NCA~\cite{https://doi.org/10.48550/arxiv.2203.14832}, that is completely algebraic.
The attractive feature of NCA is that it provides nested bases, in that the row and column bases of an admissible matrix of interaction between two boxes are synthesized from the row and column bases of their respective children.
The advantage of constructing nested bases is that it speeds up the algorithm.

In~\cite{https://doi.org/10.48550/arxiv.2203.14832}, NCA is developed for non-oscillatory kernels.
In this section, we summarise important details of~\cite{https://doi.org/10.48550/arxiv.2203.14832} and adapt it to the 2D Helmholtz kernel.
We choose to describe it with a uniform tree for pedagogical reasons. It is to be noted that it is readily extendable to a level-restricted tree.

\subsection{Preliminaries}
Let the index sets of the matrix $A$ be $I\times J$.
Matrix $A$ holds the pairwise interactions of points in $\mathcal{X}$, i.e., for $i\in I$ and $j\in J$, the $(i,j)^{th}$ entry of $A$ is the contribution of point $x_{j}$ at the point $x_{i}$, as defined in Eq.~\eqref{eq:13}.
For a box $B$, let sets $t^{B}$ and $s^{B}$ be defined as,
\begin{align}
  t^{B} &= \{i:i\in I \text{ and } x_{i}\in B\} \text{and}\\
  s^{B} &= \{j:j\in J \text{ and } x_{j}\in B\}.
\end{align}
Let the matrix sub-block representing the contribution of points $\{x_{j}\}_{j\in s^{Y}}$ at points $\{x_{i}\}_{i\in t^{X}}$ be denoted by $A_{t^{X}s^{Y}}$, i.e., the $(i,j)^{th}$ entry of $A_{t^{X}s^{Y}}$ is $A(t^{X}(i),s^{Y}(j))$.

\begin{figure}[H]
     \begin{center}
     \includegraphics[width=0.5\linewidth]{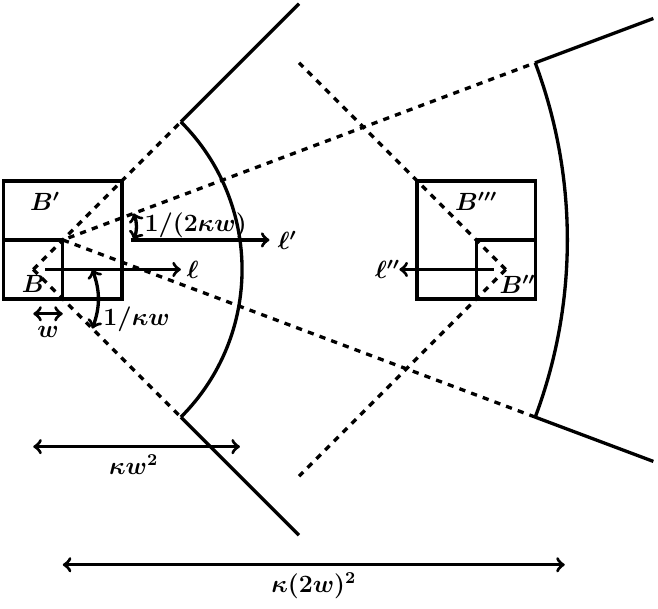}
     \caption{Illustration of box $B$ in high frequency regime, where $(B,\ell)\in\mathcal{C}(B',\ell')$ and $(B'',\ell'')\in\mathcal{IL}(B,\ell)$.}
     \label{fig:HFR_IL}
 \end{center}
 \end{figure}

We define some notations in Table~\ref{table:Notations} which will be used in the rest of the article.
In Figure~\ref{fig:HFR_IL} we illustrate a box $B$ in high frequency regime, its far-field region in direction $\ell$, its parent $B'$, and the far-field region of $B'$ in direction $\ell'$ where $(B,\ell)\in\mathcal{C}(B',\ell')$. In the same figure we also illustrate box $B''$, its parent $B'''$ and direction $\ell''$, where $(B'',\ell'')\in\mathcal{IL}(B,\ell)$. It is to be noted that $B'''\in \mathcal{N}(B')$, $B''\not\in \mathcal{N}(B)$, $B''$ falls in the cone in direction $\ell$ of $B$ and $B$ falls in the cone in direction $\ell''$ of $B''$.

\begin{table}[H]
\centering
  \begin{tabularx}{\textwidth}{|l|X|}
    \hline
    $B$ & Box at a certain level in the tree\\ \hline
    $\mathcal{C}(B)$ & $\{B':\text{ }B'\text{ is a child of }B\}$ \\ \hline
    $\mathcal{C}(B,\ell)$ & For a box $B$ whose children are in the high frequency regime, $\mathcal{C}(B,\ell) = \{(B',\ell'):B'\in C(B)$; $\ell'$ is the direction associated with the cones of its children such that the cone in direction $\ell$ of $B$ is within the cone in direction $\ell'$ of its children\}. For its illustration we refer the readers to~\cite{engquist2009fast} \\ \hline
    $\mathcal{N}(B)$ & Neighbors of $B$; For a box $B$ in the low frequency regime, it consists of boxes at the same tree level as $B$ that do not follow the admissibility condition for low-rank in low frequency regime. For a box $B$ in the high frequency regime, it consists of boxes at the same tree level as $B$ that do not follow the distance condition of the directional admissibility condition for low-rank in high frequency regime.\\ \hline
    $\mathcal{IL}(B)$ & Interaction list of a box $B$ in the low frequency regime that consists of children of $B$'s parent's neighbors that are not its neighbors.\\ \hline
    $\mathcal{IL}(B,\ell)$ & Interaction list of a box $B$ in the high frequency regime in direction $\ell\in L(B)$ that consists of pairs of boxes and directions of the form $(A,\ell')$ such that $A$ falls in the cone in direction $\ell$ of $B$ and $B$ falls in the cone in direction $\ell'$ of $A$, and $A$ is a child of a neighbor of $B$'s parent but not a neighbor of $B$.\\ \hline
  \end{tabularx}
  \caption{List of notations followed in the rest of the article}
    \label{table:Notations}
 \end{table}

 \subsection{Construction of low-rank approximations}\label{sec:constructionNCA}
Consider boxes $X$ and $Y$ in the low frequency regime, where $Y\in\mathcal{IL}(X)$.
The low-rank approximation of sub-block $A_{t^{X}s^{Y}}$, as $U_{X}S_{X,Y}V_{Y}^{*}$, using NCA takes the form~\cite{https://doi.org/10.48550/arxiv.2203.14832}:
\begin{equation} \label{eq:NCA}
  A_{t^{X}s^{Y}} \approx \underbrace{A_{t^{X}s^{X,i}} (A_{t^{X,i}s^{X,i}})^{-1}}_{U_{X}} \underbrace{A_{t^{X,i}s^{Y,o}}}_{S_{X,Y}} \underbrace{(A_{t^{Y,o}s^{Y,o}})^{-1} A_{t^{Y,o}s^{Y}}}_{V_{Y}^{*}}
\end{equation}
where $t^{X,i}\subset t^{X}$, $s^{X,i}\subset \mathcal{F}^{X,i}$, $t^{Y,o}\subset \mathcal{F}^{Y,o}$ and $s^{Y,o}\subset s^{Y}$ are termed \textit{pivots}. And $\mathcal{F}^{X,i}$ and $\mathcal{F}^{Y,o}$ are defined as
\begin{align}
  \mathcal{F}^{X,i} &= \{s^{X'}:X'\in \mathcal{IL}(X)\} \text{ and}\label{eq:Fj}\\
  \mathcal{F}^{Y,o} &= \{t^{Y'}:Y'\in \mathcal{IL}(Y)\}.\label{eq:Fi}
\end{align}
$t^{X,i}$ and $s^{X,i}$ are termed the \textit{incoming row pivots} and \textit{incoming column pivots} of $X$ respectively.
And $t^{Y,o}$ and $s^{Y,o}$ are termed the \textit{outgoing row pivots} and \textit{outgoing column pivots} of $Y$ respectively. Matrices $U_{X}$ and $V_{Y}^{*}$ are termed the column and row bases of boxes $X$ and $Y$ respectively.

Similarly, the construction of low-rank approximation of the sub-block $A_{t^{Y}s^{X}}$ involves pivots $t^{Y,i}$, $s^{Y,i}$, $t^{X,o}$ and $s^{X,o}$ and bases $U_{Y}$ and $V_{X}^{*}$. So each box $B$ in the low frequency regime is associated with pivots $t^{X,i}$, $s^{X,i}$, $t^{X,o}$ and $s^{X,o}$, column basis $U_{X}$, and row basis $V_{X}^{*}$.

$t^{X,i}$ and $s^{X,o}$ represent the points lying in the box $X$ and hence are also termed the \textit{self pivots} of $X$.
$s^{X,i}$ and $t^{X,o}$ represent the points in the far-field region of $X$ and hence are also termed the \textit{far-field pivots} of $X$.
For more details on the approximation, we refer the readers to~\cite{https://doi.org/10.48550/arxiv.2203.14832}.

The construction of low-rank approximations in the high frequency regime is similar to the low frequency regime case, except that the pivots and the bases are defined for each direction $\ell\in L(X)$ associated with a box $X$. Hence the pivots and the bases in the high frequency regime are said to be \textit{directional}.
Consider two boxes $X$ and $Y$ in the high frequency regime such that $Y\in \mathcal{IL}(X,\ell)$, for an $\ell\in L(X)$ and $X\in \mathcal{IL}(Y,\ell')$, for an $\ell'\in L(Y)$. The low-rank approximation in the high frequency regime or, what we call, the \textit{directional low-rank approximation} of matrix sub-block $A_{t^{X}s^{Y}}$, as $U_{X}^{\ell}S_{X,Y}^{\ell,\ell'}{V_{Y}^{\ell'}}^{*}$, using NCA takes the form:
\begin{equation} \label{eq:NCA_HFR}
  A_{t^{X}s^{Y}} \approx \underbrace{A_{t^{X}s^{X,i,\ell}} (A_{t^{X,i,\ell}s^{X,i,\ell}})^{-1}}_{U_{X}^{\ell}} \underbrace{A_{t^{X,i,\ell}s^{Y,o,\ell'}}}_{S_{X,Y}^{\ell,\ell'}} \underbrace{(A_{t^{Y,o,\ell'}s^{Y,o,\ell'}})^{-1} A_{t^{Y,o,\ell'}s^{Y}}}_{{V_{Y}^{\ell'}}^{*}}
\end{equation}
where $t^{X,i,\ell}\subset t^{X}$, $s^{X,i,\ell}\subset \mathcal{F}^{X,i,\ell}$, $t^{Y,o,\ell'}\subset \mathcal{F}^{Y,o,\ell'}$ and $s^{Y,o,\ell'}\subset s^{Y}$ are termed pivots. And $\mathcal{F}^{X,i,\ell}$ and $\mathcal{F}^{Y,o,\ell'}$ are defined as
\begin{align}
  \mathcal{F}^{X,i,\ell} &= \{s^{X'}:X'\in \mathcal{IL}(X,\ell)\}, \text{ and}\label{eq:Fj_HFR}\\
  \mathcal{F}^{Y,o,\ell'} &= \{t^{Y'}:Y'\in \mathcal{IL}(Y,\ell')\}.\label{eq:Fi_HFR}
\end{align}
$t^{X,i,\ell}$ and $s^{X,i,\ell}$ are termed the \textit{incoming row pivots} and \textit{incoming column pivots} of $X$ in direction $\ell$ respectively.
And we term $t^{Y,o,\ell'}$ and $s^{Y,o,\ell'}$ as the \textit{outgoing row pivots} and \textit{outgoing column pivots} of $Y$ in direction $\ell'$ respectively. Matrices $U_{X}^{\ell}$ and ${V_{Y}^{\ell'}}^{*}$ are termed the column and row bases of boxes $X$ and $Y$ in directions $\ell$ and $\ell'$ respectively.

Similarly, the construction of low-rank approximation of the sub-block $A_{t^{Y}s^{X}}$ involves pivots $t^{Y,i,\ell'}$, $s^{Y,i,\ell'}$, $t^{X,o,\ell}$ and $s^{X,o,\ell}$ and bases $U_{Y}^{\ell'}$ and ${V_{X}^{\ell}}^{*}$.
So each box and direction pair $(B,\ell\in L(B))$ in the high frequency regime, is associated with pivots $t^{X,i,\ell}$, $s^{X,i,\ell}$, $t^{X,o,\ell}$ and $s^{X,o,\ell}$, column basis $U_{X}^{\ell}$, and row basis ${V_{X}^{\ell}}^{*}$.

$t^{X,i,\ell}$ and $s^{X,o,\ell}$ represent the points lying in box $X$, that are used to construct the directional low-rank approximations of the matrix sub-blocks $A_{t^{X}s^{Y}}$ and $A_{t^{Y}s^{X}}$ and hence are also termed the \textit{self pivots} of $X$ in direction $\ell$.
$s^{X,i,\ell}$ and $t^{X,o,\ell}$ are termed the \textit{far-field pivots} of $X$ in direction $\ell$, as they represent the points that lie in the far-field region of $X$ in direction $\ell$.

\subsubsection{Construction of Nested Bases}
The column and row bases of boxes in both the high and low frequency regimes are constructed in a nested fashion: by expressing the bases of a non-leaf box in terms of the bases of its children.
We now describe the construction of bases, in four different possible scenarios.
\begin{enumerate}
  \item
\textbf{For a leaf box $B$,}
    its column and row bases, also termed the L2P and P2M translation operators of $B$ (the terminology used with FMM), are given by
    \begin{equation}
    U_{B} := A_{t^{B}s^{B,i}} A_{t^{B,i} s^{B,i}}^{-1} \hspace{5mm} \text{and} \hspace{5mm}
    V_{B}^{*} := A_{t^{B,o} s^{B,o}}^{-1} A_{t^{B,o}s^{B}}.
    \end{equation}
  \item
  \textbf{For a parent box $B$ in the low frequency regime (LFR),}
  its column and row bases are given by
  \begin{equation}
    U_{B} = \begin{bmatrix}
            U_{B_{1}} & 0 & 0 & 0\\
            0 & U_{B_{2}} & 0 & 0\\
            0 & 0 & U_{B_{3}} & 0\\
            0 & 0 & 0 & U_{B_{4}}
            \end{bmatrix}
            \begin{bmatrix}
            C_{B_{1}B}\\
            C_{B_{2}B}\\
            C_{B_{3}B}\\
            C_{B_{4}B}\\
            \end{bmatrix}
            \hspace{5mm} \text{and} \hspace{5mm}
            V_{B} = \begin{bmatrix}
                    V_{B_{1}} & 0 & 0 & 0\\
                    0 & V_{B_{2}} & 0 & 0\\
                    0 & 0 & V_{B_{3}} & 0\\
                    0 & 0 & 0 & V_{B_{4}}
                    \end{bmatrix}
                    \begin{bmatrix}
                    {T_{BB_{1}}}^{*}\\
                    {T_{BB_{2}}}^{*}\\
                    {T_{BB_{3}}}^{*}\\
                    {T_{BB_{4}}}^{*}\\
                    \end{bmatrix},
  \end{equation}
 where $\{B_{c}\}_{c=1}^{4}\in \mathcal{C}(B)$
and matrices $\{C_{B_{c}B}\}_{c=1}^{4}$ and $\{{T_{BB_{c}}}^{*}\}_{c=1}^{4}$, termed the L2L and M2M translation operators of $B$, take the following form
  \begin{equation}\label{eq:S_Non_Directional_L2L}
    C_{B_{c}B} =  A_{t^{B_{c},i}s^{B,i}} A_{t^{B,i} s^{B,i}}^{-1} \hspace{5mm} \text{and} \hspace{5mm}
    T_{BB_{c}} = A_{t^{B,o} s^{B,o}}^{-1} A_{t^{B,o}s^{B_{c},o}} \hspace{5mm} \forall c\in\{1,2,3,4\}.
  \end{equation}

  \item
  \textbf{For a parent box $B$ in the high frequency regime (HFR) and children in the LFR.}
  When the transition from the high to low frequency regime happens, a box $B$ at a parent level has bases defined for each direction $\ell\in L(B)$, whereas its children $\{B_{i}\}_{i=1}^{4}\in \mathcal{C}(B)$ have the \textit{non-directional bases}. The column and row bases of such a box $B$ in direction $\ell$ are given by
  \begin{equation}
    U_{B}^{\ell} = \begin{bmatrix}
            U_{B_{1}} & 0 & 0 & 0\\
            0 & U_{B_{2}} & 0 & 0\\
            0 & 0 & U_{B_{3}} & 0\\
            0 & 0 & 0 & U_{B_{4}}
            \end{bmatrix}
            \begin{bmatrix}
            C_{B_{1}B}^{\ell}\vspace{0.5mm}\\
            C_{B_{2}B}^{\ell}\vspace{0.5mm}\\
            C_{B_{3}B}^{\ell}\vspace{0.5mm}\\
            C_{B_{4}B}^{\ell}\vspace{0.5mm}\\
            \end{bmatrix}
            \hspace{5mm} \text{and} \hspace{5mm}
            V_{B}^{\ell} = \begin{bmatrix}
                    V_{B_{1}} & 0 & 0 & 0\\
                    0 & V_{B_{2}} & 0 & 0\\
                    0 & 0 & V_{B_{3}} & 0\\
                    0 & 0 & 0 & V_{B_{4}}
                    \end{bmatrix}
                    \begin{bmatrix}
                    {T_{BB_{1}}^{\ell}}^{*}\\
                    {T_{BB_{2}}^{\ell}}^{*}\\
                    {T_{BB_{3}}^{\ell}}^{*}\\
                    {T_{BB_{4}}^{\ell}}^{*}\\
                    \end{bmatrix},
  \end{equation}
  where $\{B_{c}\}_{c=1}^{4}\in \mathcal{C}(B)$
  and
  matrices $\{C_{B_{c}B}^{\ell}\}_{c=1}^{4}$ and $\{{T_{BB_{c}}^{\ell}}^{*}\}_{c=1}^{4}$, termed the directional L2L and directional M2M translation operators of $B$ in direction $\ell$ respectively, are defined as
  \begin{equation}
    C_{B_{c}B}^{\ell} =  A_{t^{B_{c},i}s^{B,i,\ell}} A_{t^{B,i,\ell} s^{B,i,\ell}}^{-1} \hspace{5mm} \text{and} \hspace{5mm}
    T_{BB_{c}}^{\ell} = A_{t^{B,o,\ell} s^{B,o,\ell}}^{-1} A_{t^{B,o,\ell}s^{B_{c},o}}  \hspace{5mm}\forall c\in\{1,2,3,4\}.
  \end{equation}
\item
\textbf{For a parent box $B$ in the HFR and children in the HFR,}
the column and row bases in direction $\ell\in L(B)$ are given by
\begin{equation}
  U_{B}^{\ell} = \begin{bmatrix}
          U_{B_{1}}^{\ell'} & 0 & 0 & 0\\
          0 & U_{B_{2}}^{\ell'} & 0 & 0\\
          0 & 0 & U_{B_{3}}^{\ell'} & 0\\
          0 & 0 & 0 & U_{B_{4}}^{\ell'}
          \end{bmatrix}
          \begin{bmatrix}
          C_{B_{1}B}^{\ell',\ell}\vspace{0.5mm}\\
          C_{B_{2}B}^{\ell',\ell}\vspace{0.5mm}\\
          C_{B_{3}B}^{\ell',\ell}\vspace{0.5mm}\\
          C_{B_{4}B}^{\ell',\ell}\vspace{0.5mm}\\
          \end{bmatrix}
          \hspace{5mm} \text{and} \hspace{5mm}
          V_{B}^{\ell} = \begin{bmatrix}
                  V_{B_{1}}^{\ell'} & 0 & 0 & 0\\
                  0 & V_{B_{2}}^{\ell'} & 0 & 0\\
                  0 & 0 & V_{B_{3}}^{\ell'} & 0\\
                  0 & 0 & 0 & V_{B_{4}^{\ell'}}
                  \end{bmatrix}
                  \begin{bmatrix}
                  {T_{BB_{1}}^{\ell,\ell'}}^{*}\\
                  {T_{BB_{2}}^{\ell,\ell'}}^{*}\\
                  {T_{BB_{3}}^{\ell,\ell'}}^{*}\\
                  {T_{BB_{4}}^{\ell,\ell'}}^{*}\\
                  \end{bmatrix}
\end{equation}
where $\{(B_{c},\ell')\}_{c=1}^{4}\in \mathcal{C}(B, \ell)$
 and
matrices $\{C_{B_{c}B}^{\ell',\ell}\}_{c=1}^{4}$ and $\{{T_{BB_{c}}^{\ell,\ell'}}^{*}\}_{c=1}^{4}$, termed the directional L2L and directional M2M translation operators of $B$ in direction $\ell$ respectively, are defined as
\begin{equation}
  C_{B_{c}B}^{\ell',\ell} =  A_{t^{B_{c},i,\ell'}s^{B,i,\ell}} A_{t^{B,i,\ell} s^{B,i,\ell}}^{-1} \hspace{5mm} \text{and} \hspace{5mm}
  T_{BB_{c}}^{\ell,\ell'} = A_{t^{B,o,\ell} s^{B,o,\ell}}^{-1} A_{t^{B,o,\ell}s^{B_{c},o,\ell'}}  \hspace{5mm} \forall c\in\{1,2,3,4\}.
\end{equation}

\end{enumerate}

For a box $B$ in the low frequency regime, $U_{B}$ and $V_{B}^{*}$ approximate $A_{t^{B}s^{B,i}} A_{t^{B,i} s^{B,i}}^{-1}$ and \\$A_{t^{B,o} s^{B,o}}^{-1} A_{t^{B,o}s^{B}}$ respectively.
Similarly for a box $B$ in the high frequency regime and direction $\ell\in L(B)$, $U_{B}^{\ell}$ and ${V_{B}^{\ell}}^{*}$ approximate $A_{t^{B}s^{B,i,\ell}} A_{t^{B,i,\ell} s^{B,i,\ell}}^{-1}$ and $A_{t^{B,o,\ell} s^{B,o,\ell}}^{-1} A_{t^{B,o,\ell}s^{B}}$ respectively.
We refer the readers to~\cite{bebendorf2012constructing,bebendorf2015wideband} for the error estimates.

\subsubsection{Nested Pivots for NCA}\label{section:nestedPivots}
We now describe the method to compute pivots, the key step of NCA.

\textbf{Low frequency regime.}
For a box $B$ in the low frequency regime, its self pivots $t^{B,i}$ and $s^{B,o}$ are chosen from $t^{B}$ and $s^{B}$ respectively. And its far-field pivots $s^{B,i}$ and $t^{B,o}$ are chosen from $\mathcal{F}^{B,i}$ and $\mathcal{F}^{B,o}$ respectively. Or equivalently, for a box $B$, the search space of its self pivots is itself and the search space of its far-field pivots is the union of boxes in its interaction list.

\textbf{High frequency regime.}
For a box $B$ in the high frequency regime and direction $\ell\in L(B)$, the self pivots $t^{B,i,\ell}$ and $s^{B,o,\ell}$ are chosen from $t^{B}$ and $s^{B}$ respectively. And the far-field pivots $s^{B,i,\ell}$ and $t^{B,o,\ell}$ are chosen from $\mathcal{F}^{B,i,\ell}$ and $\mathcal{F}^{B,o,\ell}$ respectively. Or equivalently, for a box $B$ and direction $\ell\in L(B)$, the search space of self pivots is itself and the search space of far-field pivots is the union of boxes in its interaction list in direction $\ell$.
To find pivots in the high frequency regime we follow the same steps as that of the low frequency regime except that the pivots are computed for all directions $\ell\in L(B)$ associated with a box $B$.

The pivots for boxes in the low and high frequency regimes are computed in a nested fashion: The pivots at a parent level of the quad-tree are computed from the pivots at its child level.
One needs to traverse the tree upwards (starting at leaf boxes) and follow the two steps described below to find the pivots by recursion.

\begin{enumerate}
  \item
  The first step in identifying nested pivots is described below for four different possible scenarios. \\
  For all leaf boxes $B$, construct sets
  \begin{equation}
    \tilde{t}^{B,i}:=t^{B}, \hspace{5mm}  \hspace{5mm} \tilde{s}^{B,i}:=\bigcup_{B'\in \mathcal{IL}(B)} s^{B'},
  \end{equation}
  \begin{equation}
     \tilde{t}^{B,o}:=\bigcup_{B'\in \mathcal{IL}(B)}t^{B'} \hspace{5mm} \text{and} \hspace{5mm}
     \tilde{s}^{B,o}:=s^{B}.
  \end{equation}
  For all non-leaf boxes $B$ in low frequency regime, construct sets
  \begin{equation}
    \tilde{t}^{B,i}:=\bigcup_{B'\in\mathcal{C}(B)} t^{B',i}, \hspace{5mm}  \hspace{5mm} \tilde{s}^{B,i}:=\bigcup_{B'\in \mathcal{IL}(B)} \hspace{1mm} \bigcup_{B''\in \mathcal{C}(B')} s^{B'',o},
  \end{equation}
  \begin{equation}
     \tilde{t}^{B,o}:=\bigcup_{B'\in \mathcal{IL}(B)} \hspace{1mm} \bigcup_{B''\in \mathcal{C}(B')} t^{B'',i} \hspace{5mm} \text{and} \hspace{5mm}
     \tilde{s}^{B,o}:=\bigcup_{B'\in\mathcal{C}(B)}s^{B',o}.
  \end{equation}
For all boxes $B$ in the high frequency regime whose children are in the low frequency regime, construct the following sets for all directions $\ell\in L(B)$
\begin{equation}
  \tilde{t}^{B,i,\ell}:=\bigcup_{B'\in\mathcal{C}(B)} t^{B',i}, \hspace{5mm}  \hspace{5mm} \tilde{s}^{B,i,\ell}:=\bigcup_{(B',\ell')\in \mathcal{IL}(B,\ell)}  \hspace{1mm} \bigcup_{B''\in \mathcal{C}(B')} s^{B'',o},
\end{equation}
\begin{equation}
   \tilde{t}^{B,o,\ell}:=\bigcup_{(B',\ell')\in \mathcal{IL}(B,\ell)}  \hspace{1mm} \bigcup_{B''\in \mathcal{C}(B')} t^{B'',i} \hspace{5mm} \text{and} \hspace{5mm}
   \tilde{s}^{B,o,\ell}:=\bigcup_{B'\in\mathcal{C}(B)}s^{B',o}.
\end{equation}

  For all boxes $B$ in the high frequency regime whose children are also in the high frequency regime, construct the following sets for all directions $\ell\in L(B)$
  \begin{equation}
    \tilde{t}^{B,i,\ell}:=\bigcup_{(B',\ell')\in\mathcal{C}(B,\ell)} t^{B',i,\ell'}, \hspace{5mm}  \hspace{5mm} \tilde{s}^{B,i,\ell}:=\bigcup_{(B',\ell')\in \mathcal{IL}(B,\ell)}  \hspace{1mm} \bigcup_{(B'',\ell'')\in \mathcal{C}(B',\ell')} s^{B'',o,\ell''},
  \end{equation}
  \begin{equation}
     \tilde{t}^{B,o,\ell}:=\bigcup_{(B',\ell')\in \mathcal{IL}(B,\ell)}  \hspace{1mm} \bigcup_{(B'',\ell'')\in \mathcal{C}(B',\ell')} t^{B'',i,\ell''} \hspace{5mm} \text{and} \hspace{5mm}
     \tilde{s}^{B,o,\ell}:=\bigcup_{(B',\ell')\in\mathcal{C}(B,\ell)}s^{B',o,\ell'}.
  \end{equation}
\item
For all boxes $B$ in the low frequency regime, perform partially pivoted ACA~\cite{rjasanow2002adaptive,zhao2005adaptive} with accuracy $\epsilon_{NCA}$ on the matrix $A_{\tilde{t}^{B,i}\tilde{s}^{B,i}}$ to find the row and column pivots, which are then assigned to pivots $t^{B,i}$ and $s^{B,i}$ respectively. Similarly perform partially pivoted ACA on the matrix $A_{\tilde{t}^{B,o}\tilde{s}^{B,o}}$ to find the pivots $t^{B,o}$ and $s^{B,o}$.
For all boxes $B$ in the high frequency regime and all directions $\ell\in L(B)$, perform partially pivoted ACA with accuracy $\epsilon_{NCA}$ on the matrix $A_{\tilde{t}^{B,i,\ell}\tilde{s}^{B,i,\ell}}$ to find the pivots $t^{B,i,\ell}$ and $s^{B,i,\ell}$. Similarly perform partially pivoted ACA on the matrix $A_{\tilde{t}^{B,o,\ell}\tilde{s}^{B,o,\ell}}$ to find the pivots $t^{B,o,\ell}$ and $s^{B,o,\ell}$.
\end{enumerate}

\begin{remark}
  The P2M/L2P translation operators of leaf boxes and the M2M/L2L translation operators of non-leaf boxes can be obtained as by-products of the pivot-choosing routine and need not be computed separately. For more details, we refer the readers to~\cite{https://doi.org/10.48550/arxiv.2203.14832}.
\end{remark}

\section{Directional Algebraic FMM (DAFMM)}\label{section:DFMM}
In this section, we develop the DAFMM that is based on NCA, an efficient algorithm to compute matrix-vector products involving the 2D Helmholtz kernel.
We choose to describe it with a uniform tree for pedagogical reasons. It is to be noted that it is readily extendable to a level-restricted tree.

Let the vector to be applied to the matrix be $v$.
Let the result of the matrix-vector product be $u$.
For a box $B$, let $v^{B}$ and $u^{B}$ be the sliced $v$ and $u$ vectors that represent the weights on its grid points respectively. In the steps below, that describe the algorithm of DAFMM, we follow the usual FMM terminology~\cite{GREENGARD1987325,greengard1988rapid,fong2009black} coupled with its \textit{directional} counterparts.

\begin{enumerate}
  \item
  Construct a quad-tree over the computational domain.
  And in the high frequency regime, construct a hierarchy of cones that subdivides the far-field region of a box into conical regions as described in Section~\ref{Construction_of_cones}.
  \item
  Traverse up the tree (starting at leaf nodes) to compute pivots of all boxes, P2M and L2P operators of leaf boxes, and, M2M and L2L translation operators of non-leaf boxes in the low frequency regime.
  In the high frequency regime compute pivots and the M2M, L2L translation operators of all box and direction pairs using the steps detailed in Section~\ref{section:nestedPivots}.
  \item
  \textbf{Upward Pass: }Traverse up the tree starting at the leaf level to compute the P2M/\\M2M operation, until the level where the box and direction pairs have non-empty interaction list sets is reached.

  \textbf{Non-Directional P2M/M2M.}
  For each leaf box $B$ compute multipoles,
  \begin{equation*}
    v^{B,o} = V_{B}^{*} v^{B}.
  \end{equation*}

  For each non-leaf box $B$ in the low frequency regime, compute the multipoles by recursion,
  \begin{equation*}
    v^{B,o} = \sum_{B'\in \mathcal{C}(B)} T_{BB'} v^{B',o}.
  \end{equation*}
  \textbf{Directional M2M.}
  For each box $B$ in the high frequency regime, whose children are in the low frequency regime, iterate over each direction $\ell\in L(B)$ to compute the multipoles of $B$ in direction $\ell$ by recursion,
  \begin{equation*}
    v^{B,o,\ell} = \sum_{B'\in \mathcal{C}(B)} T^{\ell}_{BB'}v^{B',o}.
  \end{equation*}
  For each box $B$ in the high frequency regime, whose children are also in the high frequency regime, iterate over each direction $\ell\in L(B)$ to compute multipoles of $B$ in direction $\ell$ by recursion,
    \begin{equation*}
      v^{B,o,\ell} = \sum_{(B',\ell')\in \mathcal{C}(B,\ell)} T_{BB'}^{\ell,\ell'} v^{B',o,\ell'}.
    \end{equation*}

  \item
  \textbf{Transverse Pass:} consists of computing the M2L operations for all boxes at all levels.

  \textbf{Non-Directional M2L.}
  For each box $B$ in the low frequency regime compute
  \begin{equation*}
    u^{B,i} = \sum_{B'\in\mathcal{IL}(B)}A_{t^{B,i}s^{B',o}} v^{B',o}.
  \end{equation*}
  \textbf{Directional M2L.}
  For each box $B$ in the high frequency regime iterate over each direction $\ell\in L(B)$ to compute the partial local expansions,
  \begin{equation*}
    u^{B,i,\ell} = \sum_{(B',\ell')\in\mathcal{IL}(B,\ell)} A_{t^{B,i,\ell}s^{B',o,\ell'}} v^{B',o,\ell'}.
  \end{equation*}
  \item
  \textbf{Downward Pass: }Traverse down the tree starting at the level where the box and direction pairs have non-empty interaction list sets, until the leaf level is reached.

  \textbf{Directional L2L.}
  For each box $B'$ in the high frequency regime, iterate over each direction $\ell'\in L(B')$ to add the L2L computation to the local expansions by recursion,

    \begin{equation*}
      u^{B',i,\ell'} := u^{B',i,\ell'} + \sum_{\ell\in D} C^{\ell',\ell}_{B'B} u^{B,i,\ell}
    \end{equation*}
    where $(B',\ell')\in\mathcal{C}(B,\ell)$ and $D=\{\ell:(B',\ell')\in\mathcal{C}(B,\ell)\}$.

  \textbf{Non-Directional L2L/L2P.}
  For each non-leaf box $B'$ in the low frequency regime, whose parent $B$ is in the high frequency regime, add the L2L computation to the local expansions by recursion,
  \begin{equation*}
    u^{B',i} := u^{B',i} + \sum_{\ell\in L(B)} C^{\ell}_{B'B} u^{B,i,\ell}.
  \end{equation*}
  For each non-leaf box $B'$ in the low frequency regime, whose parent $B$ is also in the low frequency regime, add the L2L computation to the local expansions by recursion,
  \begin{equation*}
    u^{B',i} := u^{B',i} + C_{B'B} u^{B,i}.
  \end{equation*}
  For each leaf box $B$, perform the L2P computation to find the partial particle expansion,
  \begin{equation*}
    u^{B}_{calc} := U_{B} u^{B,i}.
  \end{equation*}

  \item
  Compute the \textbf{Near field} for each leaf box and add it to the particle expansion,
  \begin{equation*}
    u^{B} := u^{B} + \sum_{B'\in\mathcal{N}(B)}A_{t^{B}s^{B'}} v^{B'}.
  \end{equation*}
\end{enumerate}

\subsection{Complexity}
For a square domain with width $W$, $O(W^{2}) = O(N)$ and the number of levels in the quad tree is $\mathcal{O}(\log W)$.
We now state the complexities for the computations in the high frequency regime.
 \subsubsection{Time Complexity}
\begin{itemize}
\item
\textit{Finding pivots}. The cost to find pivots of a box $B$ with width $2w$ in direction $\ell\in L(B)$ in the high frequency regime is $O(1)$. The number of directions associated with a box of width $2w$ at level $\nu$ is $\mathcal{O}(w)=\mathcal{O}(\frac{W}{2^{\nu}})$ and the number of boxes at level $\nu$ is $4^{\nu}$. So the cost to find pivots of all boxes in the high frequency regime is $\sum_{\nu=0}^{\log W}4^{\nu}\mathcal{O}(\frac{W}{2^{\nu}})=\mathcal{O}(N)$.
\item
\textit{Directional M2M/L2L}: The cost to apply the M2M/L2L operator of a box $B$ with width $2w$ in direction $\ell\in L(B)$ at level $\nu$ is $\mathcal{O}(1)$. On similar lines of evaluating the cost for finding pivots, the cost to apply the directional M2M/L2L of all boxes in the high frequency regime is $\mathcal{O}(N)$.
\item
\textit{Directional M2L}: The number of box and direction pairs in the interaction list of a box $B$ with width $2w$ in direction $\ell\in L(B)$ at level $\nu$ is $O(w)$. For each element of the interaction list, the cost to apply the M2L operator is $O(1)$. On similar lines of evaluating the cost for finding pivots, the cost to apply the directional M2L of all boxes in the high frequency regime is $\sum_{\nu=0}^{\log W}4^{\nu}\mathcal{O}(\frac{W}{2^{\nu}})\mathcal{O}(w)=\mathcal{O}(N\log N)$.
\end{itemize}
The time complexity for the non-directional computations of the algorithm can be computed on similar lines as above and is equal to $\mathcal{O}(N)$. Hence the overall complexity of the algorithm is $\mathcal{O}(N\log N)$.

  \subsubsection{Memory Complexity}
  We now state the memory complexities of the pivots, directional M2M/L2L, and directional M2L for a single box and direction pair, $(B,\ell\in(B))$, in high frequency regime.
  \begin{itemize}
  \item
  \textit{Pivots}. The memory needed to store the pivots of a box $B$ in direction $\ell$ is $\mathcal{O}(1)$.
  \item
  \textit{Directional M2M/L2L}: The memory needed to store the directional M2M/L2L operator of a box $B$ in direction $\ell\in L(B)$ is $\mathcal{O}(1)$.
  \item
  \textit{Directional M2L}: For each element of the interaction list of a box $B$ in direction $\ell\in L(B)$, the cost to store the directional M2L operator is $\mathcal{O}(1)$.
  \end{itemize}
  The memory complexities of the pivots, directional M2M/L2L, and directional M2L for all boxes and their associated directions at all levels in the high frequency regime can be found on similar lines of evaluating their respective time complexities, and are equal to $\mathcal{O}(N)$, $\mathcal{O}(N)$, and $\mathcal{O}(N\log N)$ respectively.
  The memory complexity for the non-directional computations of the algorithm can be computed on similar lines as above and is equal to $\mathcal{O}(N)$. Hence the overall memory complexity of the algorithm is $\mathcal{O}(N\log N)$.

\section{Numerical Results and Discussion}\label{section:numericalResults}
We solve the system of Eq.~\eqref{eq:12a} for a wide range of scatterers. We compare the performance of three different solvers:
\begin{enumerate}
	\item
	\textit{GMRES accelerated by DAFMM without preconditioner}: We choose GMRES~\cite{saad1986gmres,saad2003iterative} to be our iterative solver since it is applicable without many constraints.
  Each iteration of GMRES involves computing a matrix-vector product. We employ DAFMM to compute these matrix-vector products.
  We shortly refer to it as the \textit{GMRES solver}.
	\item
	\textit{HODLR direct solver}: We also solve the arising linear system using the HODLR direct solver~\cite{ambikasaran2013mathcal},~\cite{Ambikasaran2019}. The goal is to compare the fast direct solver with its fast iterative counterpart.
	\item
	\textit{GMRES accelerated by DAFMM with HODLR as preconditioner}:
	To improve the convergence of our iterative solver, we use a HODLR~\cite{ambikasaran2013mathcal} based preconditioner. The HODLR solver can be used as a preconditioner by apriori fixing the rank of approximation of the off-diagonal sub-blocks. Let $\tilde{A}$ be such a HODLR approximate of $A$. Then the preconditioned system is
	\begin{equation} \label{eq:hybrid}
	 \tilde{A}^{-1}A\vec{\psi} = \tilde{A}^{-1}\vec{f}.
	\end{equation}
  We term the GMRES solver with HODLR based preconditioner as \textit{Hybrid solver}.
\end{enumerate}

To find $u^{scat}$, we discretize Eq.~\eqref{eq:5} using the same grid that we used to find $\psi$. We then use DAFMM to find $u^{scat}$ from the discretized $\psi$.

\subsection{Time and memory complexities}
We now state the time and memory complexities of the three solvers. A summary of the same is given in Table~\ref{table:complexities}.
\begin{enumerate}
  \item
  \textit{GMRES solver.}
  Since the time complexity of DAFMM is $\mathcal{O}(N\log N)$, the time complexity of our iterative solver with DAFMM and no preconditioner is $\mathcal{O}(m_{1}N\log N)$, where $m_{1}$ is the number of iterations it takes for GMRES to converge to a given accuracy.
  As the memory complexity of DAFMM is $\mathcal{O}(N\log N)$, the memory complexity of our iterative solver with DAFMM and no preconditioner is $\mathcal{O}(N\log N)$.
  \item
  \textit{HODLR direct solver.}
  The time complexities of HODLR factorization and solve are $\mathcal{O}(r^{2}N\log^{2} N)$ and $\mathcal{O}(rN\log N)$ respectively, where $r$ is the off-diagonal block rank~\cite{ambikasaran2013mathcal}.
  The memory complexity of both the HODLR factorization and solve is $\mathcal{O}(rN\log N)$. Since $r$ scales with $\sqrt{N}$ in $2$D, i) the time complexity of HODLR factorization and solve are $\mathcal{O}(N^{2}\log^{2} N)$ and $\mathcal{O}(N^{1.5}\log N)$ respectively ii) the memory complexity of both the HODLR factorization and solve is $\mathcal{O}(N^{1.5}\log N)$.
  \item
  \textit{Hybrid solver.}
  The time complexity of the Hybrid solver includes i) the time complexity of factorizing the HODLR preconditioner, $\mathcal{O}(r^{2}N\log^{2} N)$, where $r$ is fixed apriori to a small value ii) the time complexity of the solve part that includes the complexities of applying the HODLR preconditioner and DAFMM, which are\\
  $\mathcal{O}(m_{2}rN\log N)$ and $\mathcal{O}(m_{2}N\log N)$ respectively, where $m_{2}$ is the number of iterations it takes for the GMRES solver with preconditioner to converge to a given accuracy.
  The memory complexity of the Hybrid solver includes i) $\mathcal{O}(rN\log N)$ for factorizing the HODLR preconditioner ii) $\mathcal{O}(rN\log N)$ for applying the HODLR preconditioner and $\mathcal{O}(N\log N)$ for DAFMM, where $r$ is fixed apriori to a small value.
\end{enumerate}

\begin{table}[H]
  \centering
  \setlength\tabcolsep{4pt}
  \begin{tabularx}{\textwidth}{|c|c|c|c|X|c|}
    \hline
    & GMRES & \multicolumn{2}{|c|}{HODLR direct solver} & \multicolumn{2}{|c|}{Hybrid}\\ \cline{2-6}
    & Solve & Factorization & Solve & Factorization (HODLR) & Solve  \\ \hline\hline
    Time & $O(m_{1}N\log N)$ & $O(N^{2}\log^{2} N)$ & $O(N^{1.5}\log N)$ & $O(r^{2}N\log^{2} N)$ & $O(m_{2}rN\log N)$ \\ \hline
    Memory & $O(N\log N)$ & $O(N^{1.5}\log N)$ & $O(N^{1.5}\log N)$ & $O(rN\log N)$ & $O(rN\log N)$ \\ \hline
  \end{tabularx}
    \caption{Summary of the time and memory complexities of the three solvers.}
    \label{table:complexities}
 \end{table}

\subsection{Experiments}
We present a comprehensive set of validations and numerical benchmarks for the proposed algorithm. We consider a total of nine different experiments.
In experiment $1$, we validate DAFMM by solving an $N$-body problem with the Green's function of the Helmholtz equation as the kernel function.
In experiments $2-9$, we solve the system of Eq.~\eqref{eq:12a} for a wide range of scatterers.
To demonstrate the applicability of the solvers for a wide range of wavenumbers, we vary $\kappa$ from $40$ to $300$.
We use a plane wave in $x$ direction, $\exp(i\kappa x)$, as the incident field. In experiments $2-4$, we validate the DAFMM based iterative solver and illustrate its convergence. We also present the applicability of HODLR as a preconditioner. In experiment $5$, we use HODLR direct solver to solve for the scattered field.
In experiments $6-9$, we compare the CPU times of three different solvers.
All computations were carried out on a quad-core, 2.3 GHz Intel Core i5 processor with 8GB RAM.

The notations described in Table~\ref{table:Notations2} will be used in the rest of the section.

\begin{table}[H]
    \centering
  \begin{tabularx}{\textwidth}{|l|X|}
    \hline
 $N$ & The number of discretization points \\
\hline
$\hat{\psi}$ & Approximate of $\vec{\psi}$ that is computed by the iterative solver \\
 \hline
 $\epsilon_{NCA}$ & Compression tolerance of NCA; For a given $\epsilon_{NCA}$, we use the stopping criterion of the partially pivoted ACA stated in~\cite{rjasanow2002adaptive} to terminate the ACA's in the routine to choose pivots\\
 \hline
 $\epsilon_{grid}$ & Tolerance used in the adaptive discretization of grid\\
 \hline
 $\epsilon_{GMRES}$ & relative residual $||A\hat{\psi}-\vec{f}||_{2}/||\vec{f}||_{2}$ that is used as the stopping criterion for GMRES\\
 \hline
 T\textsubscript{Hf} & Time taken to factorize the matrix by HODLR direct solver\\
 \hline
 T\textsubscript{Hs} & Time taken to solve $\psi$ by HODLR direct solver\\
 \hline
 T\textsubscript{HODLR} & T\textsubscript{Hf} $+$ T\textsubscript{Hs}\\
 \hline
 T\textsubscript{GMRES} & Time taken to solve $\psi$ by GMRES solver with no preconditioner, wherein each iteration of GMRES involves applying DAFMM.\\
 \hline
 T\textsubscript{Pf} & Time taken to build (factorize) the preconditoner\\
 \hline
 T\textsubscript{Ps} & Time taken to solve $\psi$ by the Hybrid solver or the GMRES solver with preconditioner, wherein each iteration of GMRES involves applying the preconditioner and DAFMM.\\
 \hline
 T\textsubscript{Hybrid} & T\textsubscript{Pf} $+$ T\textsubscript{Ps}\\
 \hline
\end{tabularx}
\caption{List of notations followed in the rest of the section}
   \label{table:Notations2}
\end{table}

To demonstrate the accuracy of the solvers, we define an error function $E(x)$ in terms of the residual of Eq.~\eqref{eq:11c}, as
\begin{equation}\label{eq:error}
  E(x)= \frac{1}{\kappa^{2}}\left|\psi(x)+\kappa^{2}q(x)
  \sum_{\mathcal{L}}\int_{B}G_{\kappa}(x,y)\sum_{l=1}^{N_{p}}  Q^{\dagger}(l,:)\vec{\psi}^{B}b_{l}\left(\frac{x_{1}-\alpha_{1}^{B}}{\beta^{B}}, \frac{x_{2}-\alpha_{2}^{B}}{\beta^{B}}\right) dy - f(x))\right|, \hspace{2mm} x\in\Omega.
\end{equation}
The $\vec{\psi}(x)$ obtained from the solver is used to compute the approximate $\psi(x),\text{ } \forall x\in\Omega$ using Eq.~\eqref{eq:7}.
We use this approximate $\psi(x)$ and $\vec{\psi}(x)$ in Eq.~\eqref{eq:error}, to get the residual error $E(x)$.
It is to be noted that $E(x)$ is used to compute the error not just at the grid points, but $\forall \text{ }x\in\Omega$.
Further, note that the error, $E(x)$, is the \enquote{true} error in the sense, it captures the error
i) due to solver ii) due to discretization - this is because the computation of $f(x)$ in Eq.~\eqref{eq:error} is exact (upto roundoff) and \emph{does not depend on the grid}.

\subsubsection{Experiment 1: Validation of DAFMM}
\begin{figure}[H]
  \begin{center}
    \begin{subfigure}[b]{0.5\textwidth}
        \includegraphics[width=\linewidth]{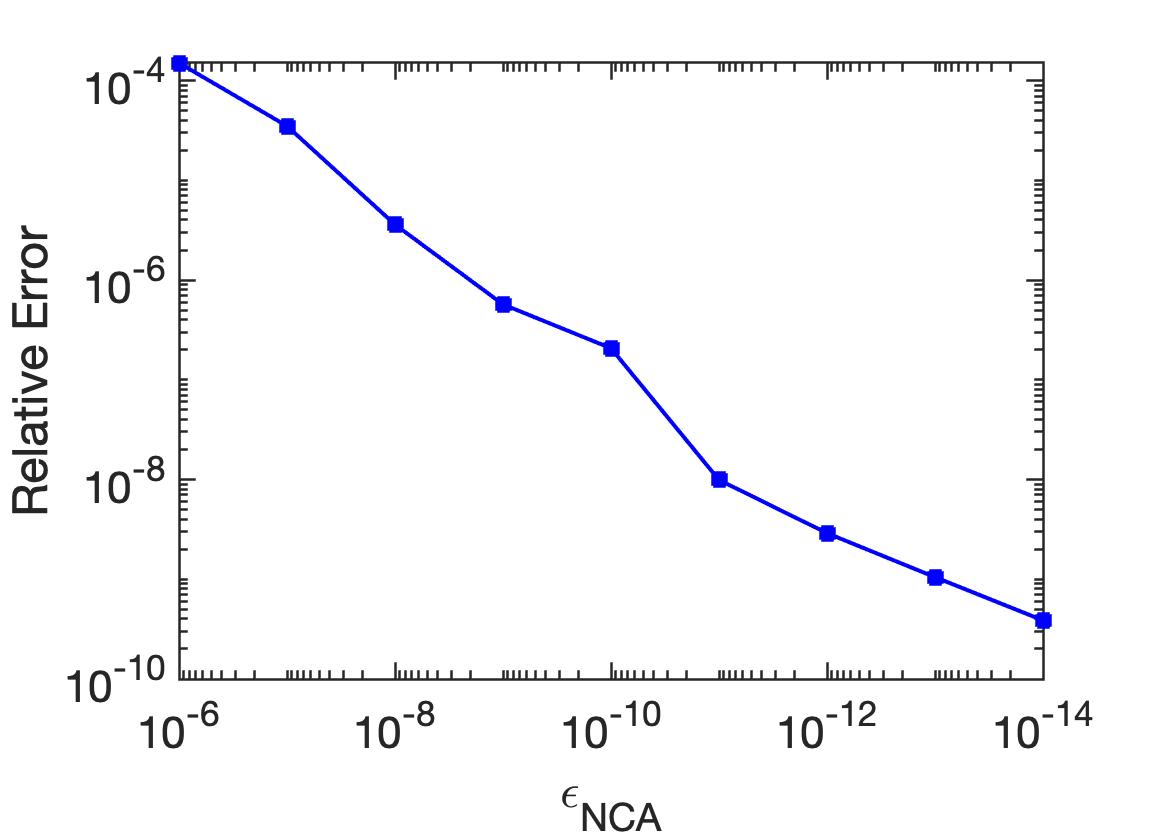}
        \caption{}
        \label{fig:errorVsEpsilon}
    \end{subfigure}%
    \begin{subfigure}[b]{0.5\textwidth}
        \includegraphics[width=\linewidth]{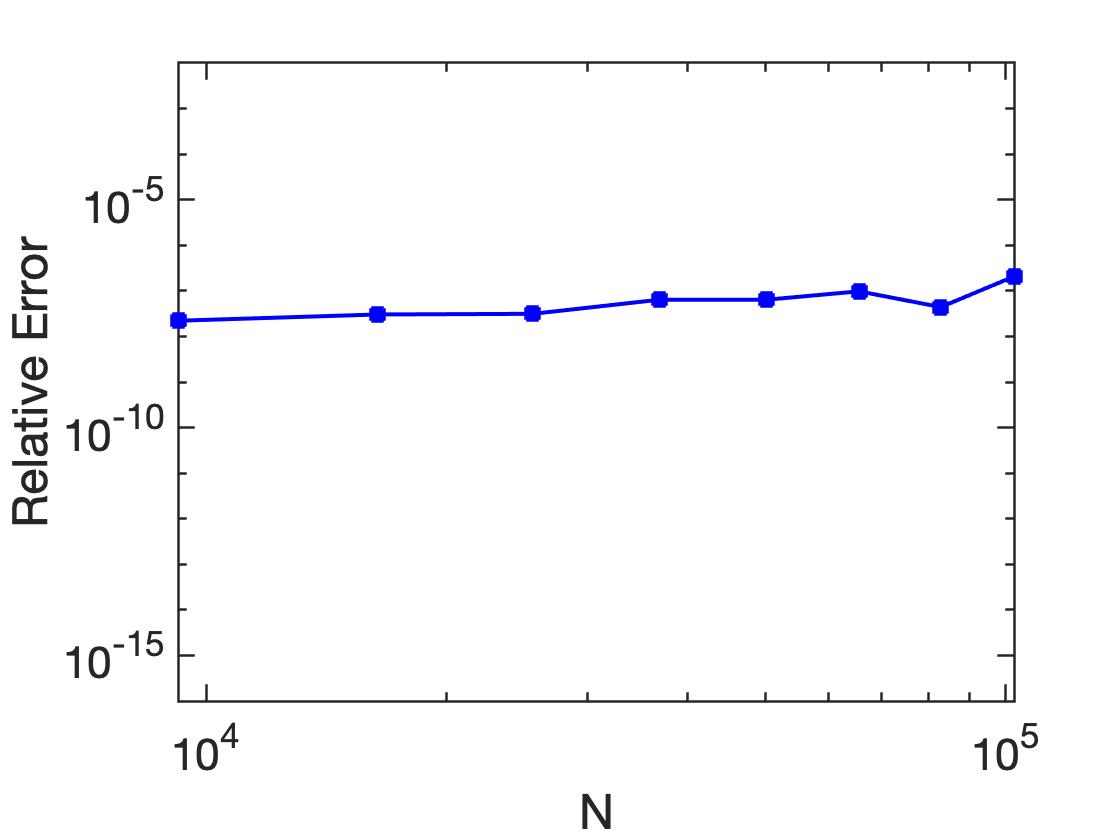}
        \caption{}
        \label{fig:errorVsN}
    \end{subfigure}%
    \caption{Results obtained with experiment 1; Relative error for the matrix-vector product a) versus $\epsilon_{NCA}$ for a system of size $102400$ b) versus $N$ keeping $\epsilon_{NCA}$ as $10^{-10}$.}
    \label{fig:DAFMMValidation}
\end{center}
\end{figure}

\begin{figure}[H]
  \begin{center}
    \begin{subfigure}[b]{0.5\textwidth}
        \includegraphics[width=\linewidth]{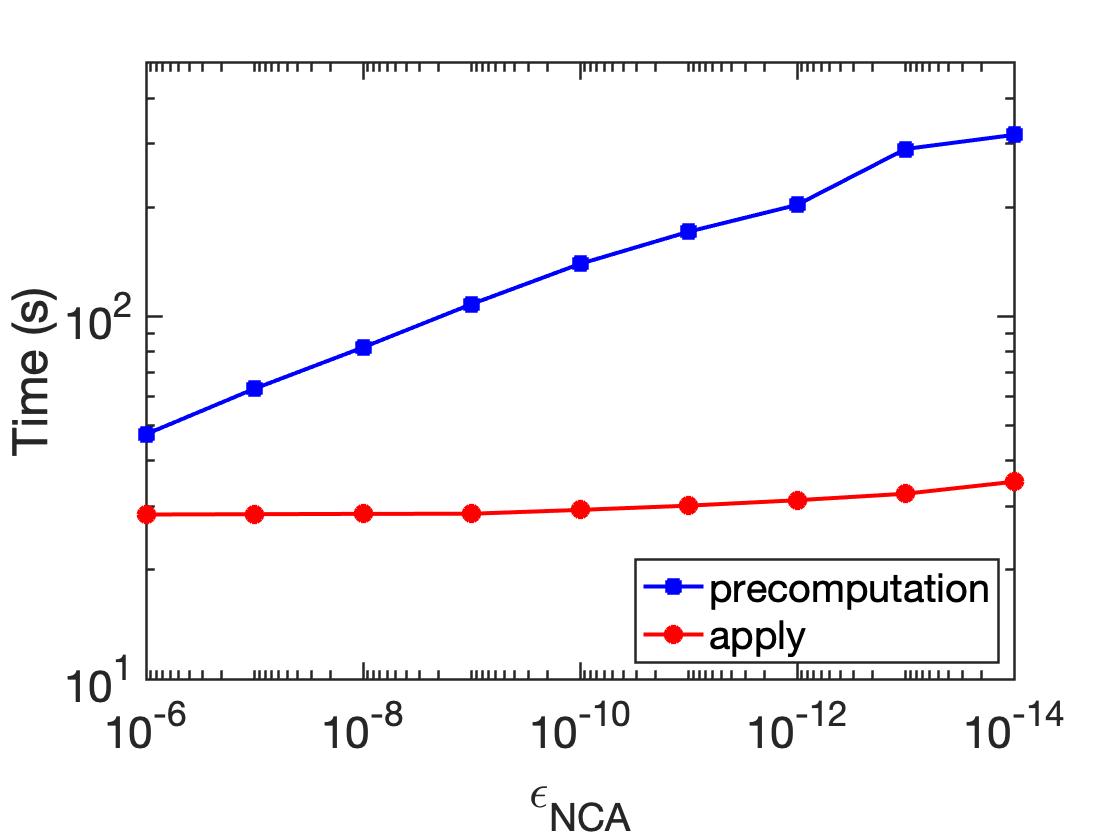}
        \caption{}
        \label{fig:timeVsEpsilon}
    \end{subfigure}%
    \begin{subfigure}[b]{0.5\textwidth}
        \includegraphics[width=\linewidth]{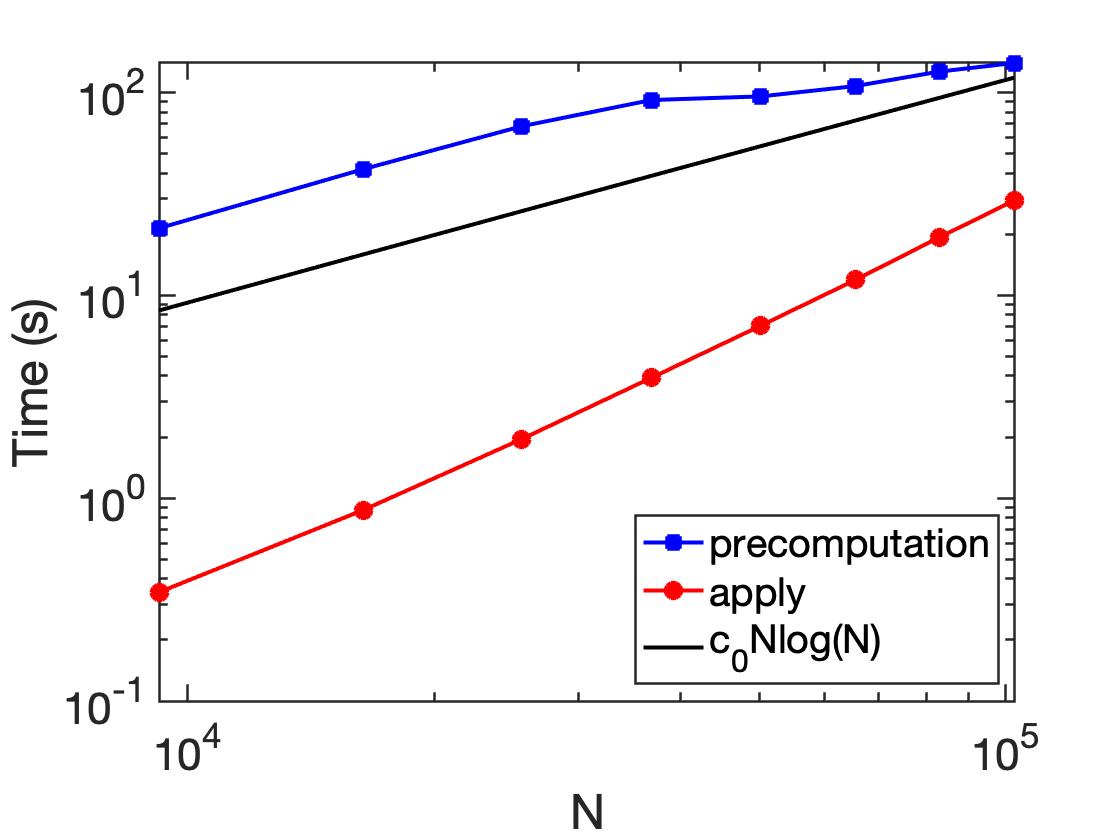}
        \caption{}
        \label{fig:timeVsN}
    \end{subfigure}%
    \caption{Results obtained with experiment 1; Precomputation time and the time taken to apply DAFMM to vector $f$ a) versus $\epsilon_{NCA}$ for a system of size $102400$ b) versus $N$ keeping $\epsilon_{NCA}$ as $10^{-10}$.}
    \label{fig:DAFMMValidation2}
\end{center}
\end{figure}

To illustrate the convergence of DAFMM we solve an $N$-body problem with the Green's function of the 2D Helmholtz equation as the kernel function. Let $f$ be a vector of $N$ charges located at points $\{p_{j}\}_{i=1}^{N}$. We compute potential $u$ at $\{p_{i}\}_{i=1}^{N}$ defined by
\begin{equation}
  u(i) = \sum_{j=1}^{N} G_{\kappa}(p_{i},p_{j})f(j)
\end{equation}
where $G_{\kappa}(x,y)=\frac{i}{4}H_{0}^{(1)}(\kappa|x-y|)$. We use the following setting for the experiment. Consider a $[-1,1]^{2}$ square domain, with $\kappa=50.0$. A uniform quad-tree is constructed such that the leaves of the tree are in the low frequency regime. In each leaf a tensor product Chebyshev grid of size $p\times p$ is considered. These grid points serve as the location of charges.
We set $p$ to $10$ and $f$ to a random vector. With these input settings, the generated system is of size $N=102400$.
We compute $\hat{u}$, an approximate of $u$, using DAFMM. The relative error $||u-\hat{u}||_{2}/||u||_{2}$ is plotted as a function of $\epsilon_{NCA}$ in Figure~\ref{fig:errorVsEpsilon}.
The precomputation time (the time taken to find pivots, form the M2M, M2L, and L2L operators in both high and low frequency regimes excluding the time taken to get the matrix entries) and the time taken to apply DAFMM to vector $f$ are plotted as a function of $\epsilon_{NCA}$ in Figure~\ref{fig:timeVsEpsilon}.

For the plots of relative error, precomputation time and apply time versus $N$, illustrated in Figures~\ref{fig:errorVsN} and~\ref{fig:timeVsN}, we vary $p$ to generate different systems sizes and keep $\epsilon_{NCA}$ constant at $10^{-10}$.

\subsubsection{Experiment 2: DAFMM accelerated GMRES solver for Gaussian contrast with $\kappa=40$}
\begin{figure}[H]
  \begin{center}
    \begin{subfigure}[b]{0.33\textwidth}
             \includegraphics[width=\linewidth]{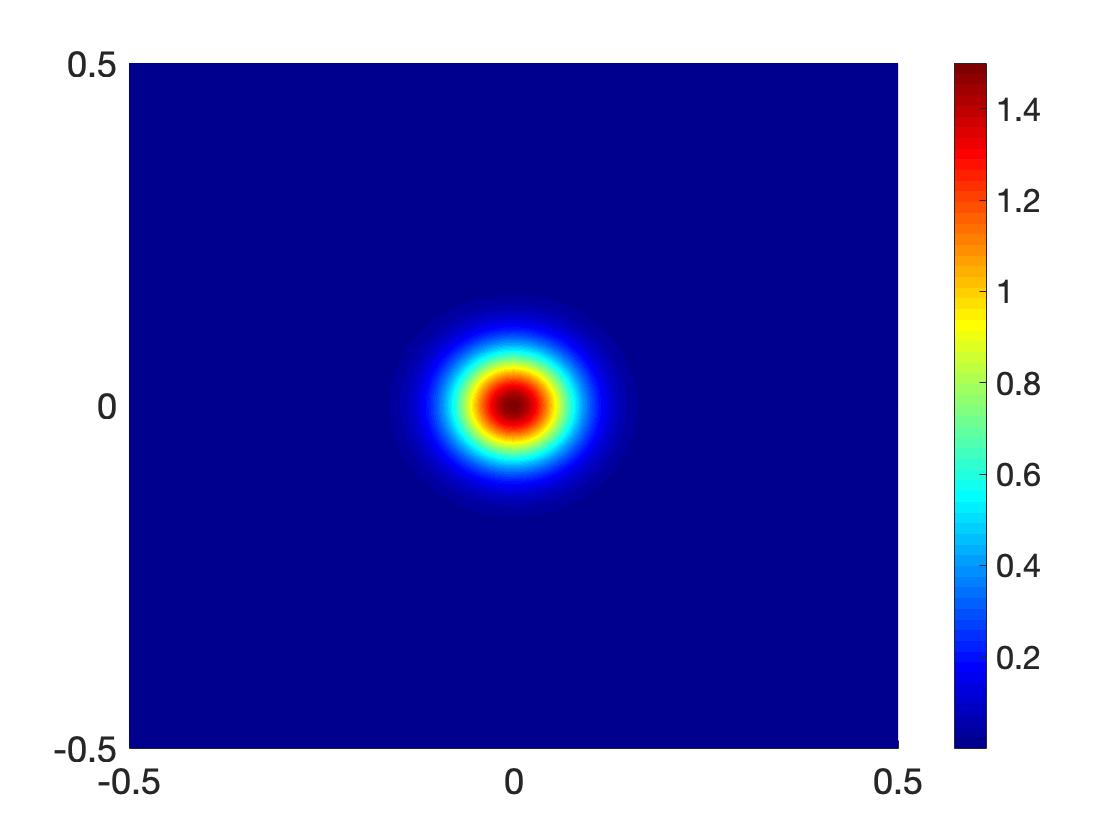}
             \caption{}
             \label{fig:GaussianContrast}
         \end{subfigure}%
         \begin{subfigure}[b]{0.33\textwidth}
             \includegraphics[width=\linewidth]{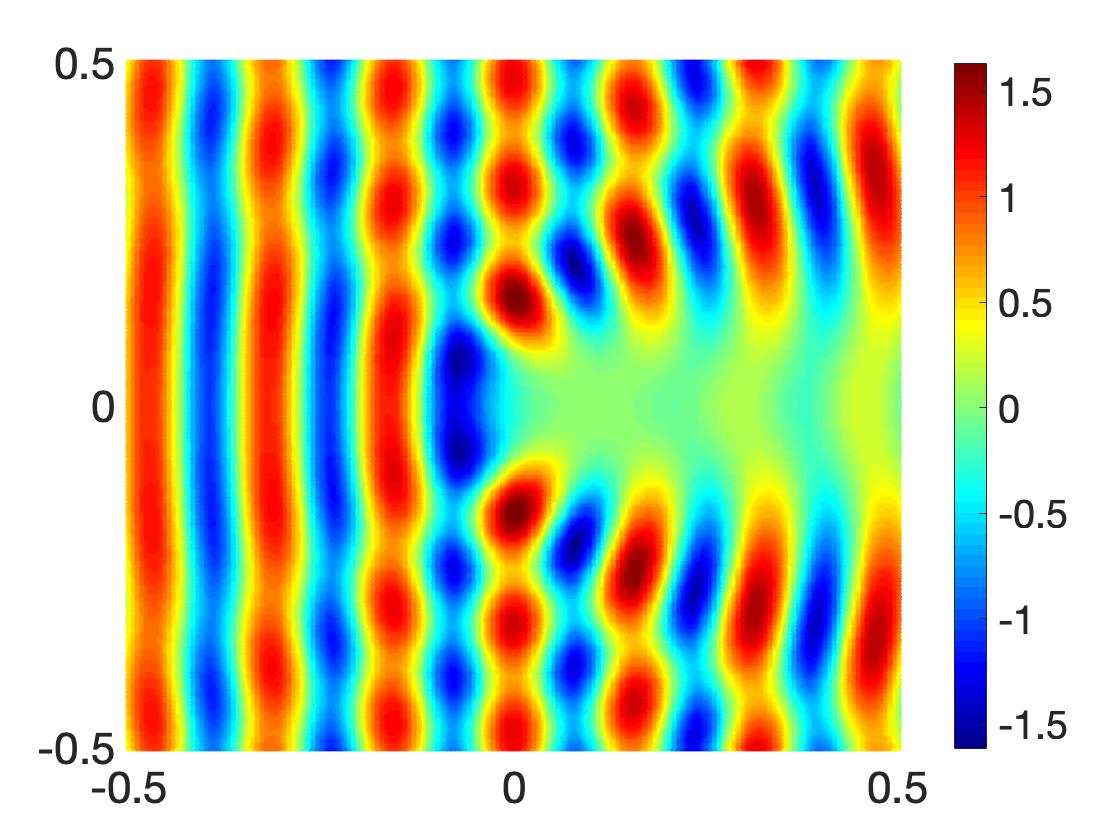}
             \caption{}
             \label{fig:GMRES40Field}
         \end{subfigure}%
         \begin{subfigure}[b]{0.33\textwidth}
             \includegraphics[width=\linewidth]{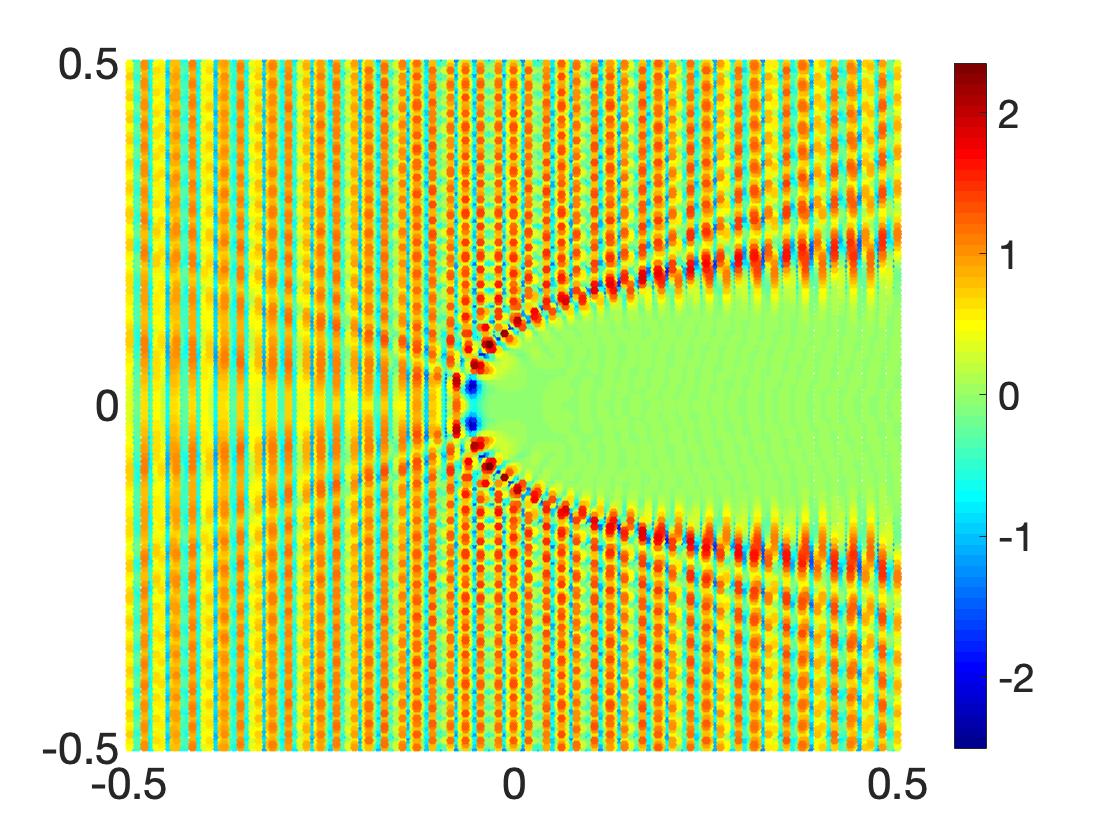}
             \caption{}
             \label{fig:GMRES300Field}
         \end{subfigure}%
         \caption{(a)Gaussian contrast; (b)Real part of the total field obtained using GMRES solver with $\kappa=40$; (c)Real part of the total field obtained using Hybrid solver with $\kappa=300$.}
         \label{fig:GaussianContrastTree}
\end{center}
\end{figure}

\begin{figure}[H]
  \begin{center}
    \includegraphics[width=0.5\linewidth]{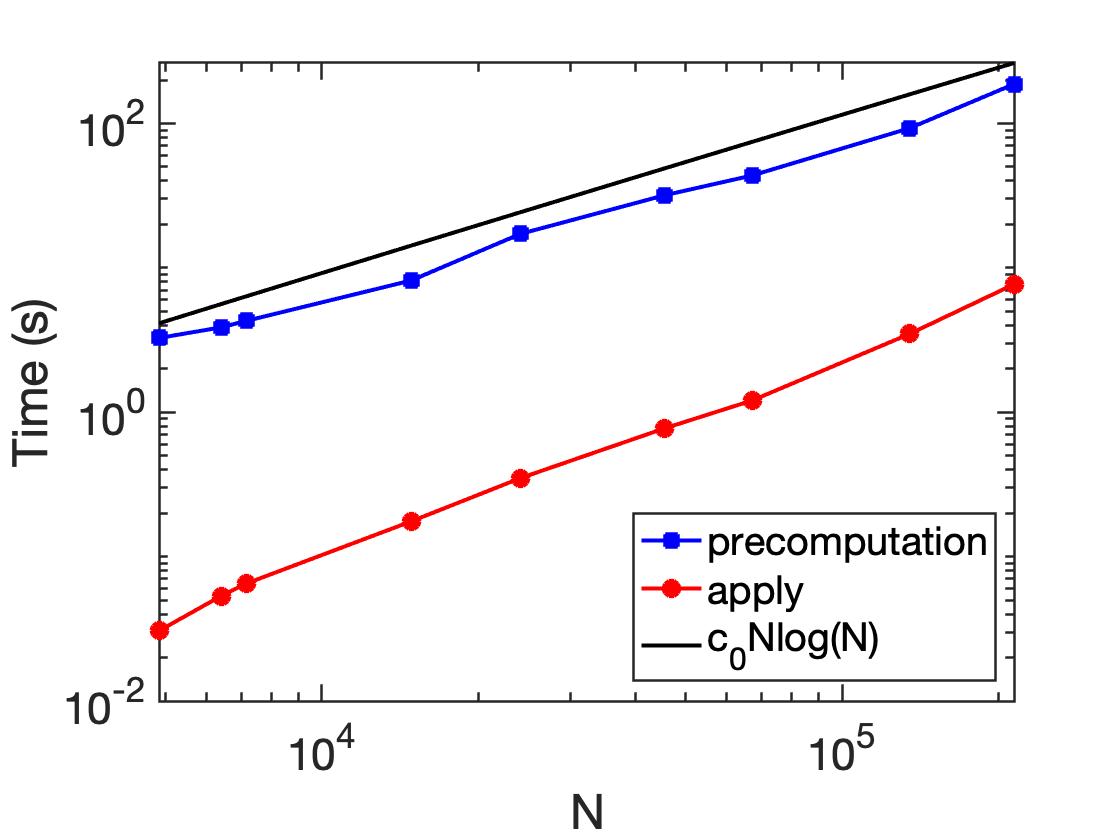}
    \caption{Results obtained with experiment 2; Precomputation time of DAFMM and the time taken to apply DAFMM to a single vector versus $N$.}
    \label{DAFMM_single}
\end{center}
\end{figure}

\begin{figure}[H]
  \begin{center}
    \begin{subfigure}[b]{0.4\textwidth}
        \includegraphics[width=\linewidth]{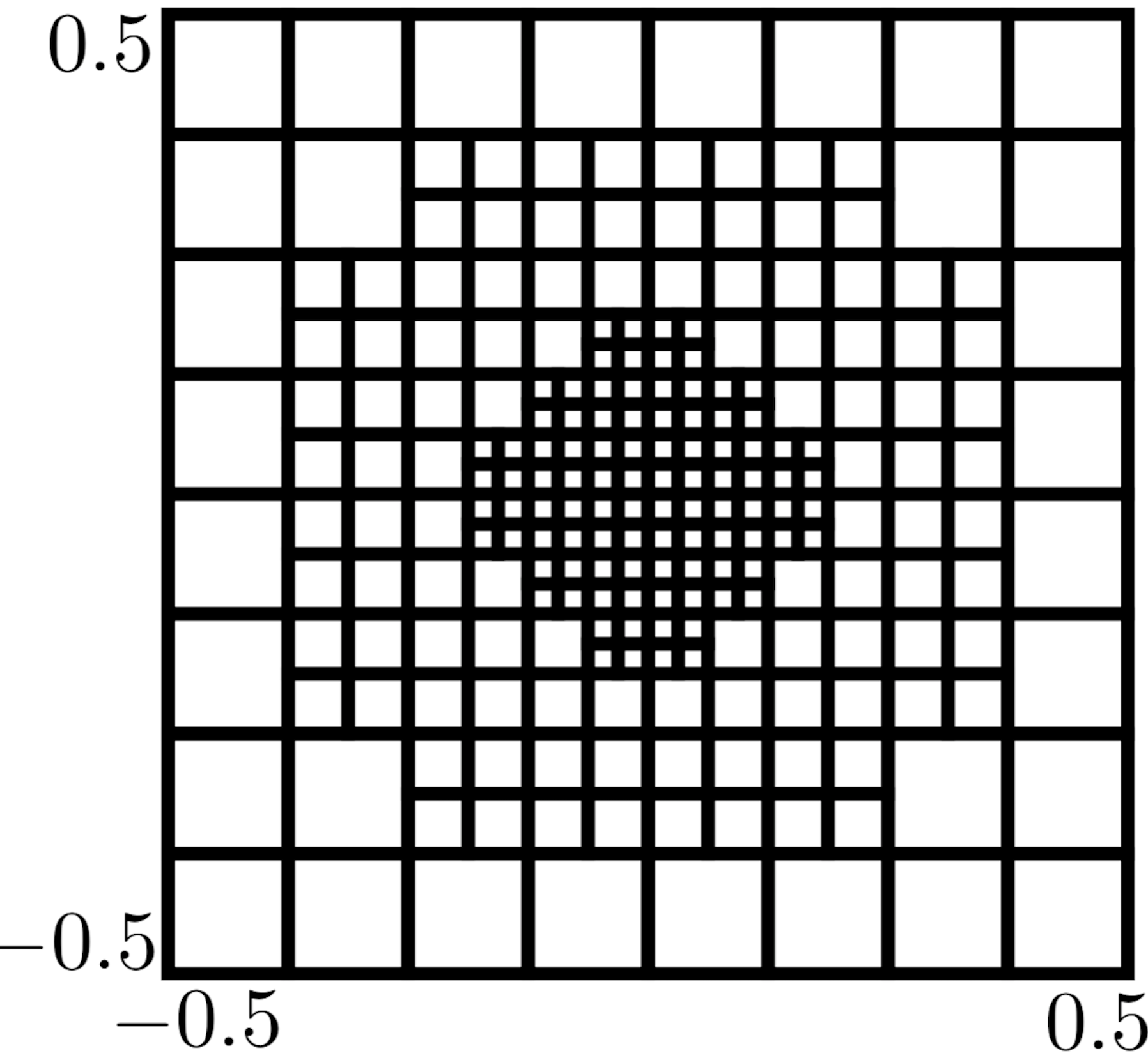}
        \caption{}
        \label{fig:grid_8}
    \end{subfigure}%
    \hfill
    \begin{subfigure}[b]{0.5\textwidth}
      \includegraphics[width=\linewidth]{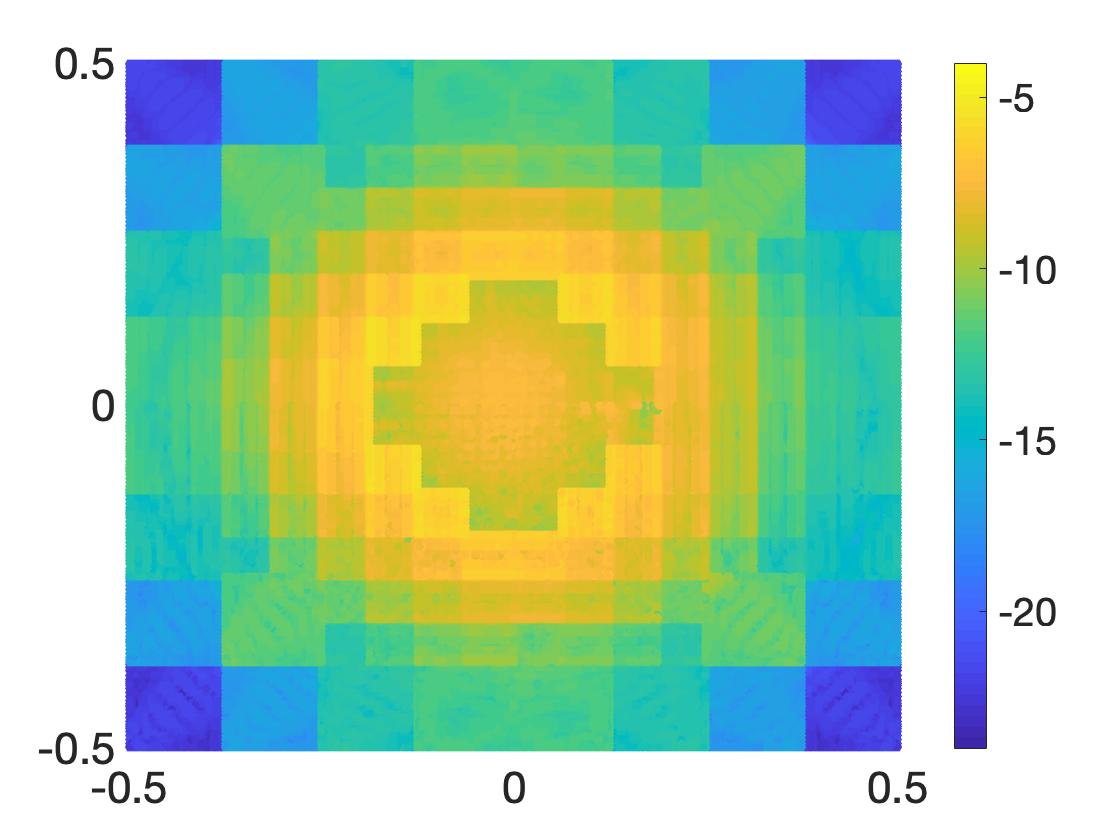}
        \caption{}
        \label{fig:error_grid_8}
    \end{subfigure}%
    \vskip\baselineskip
    \begin{subfigure}[b]{0.4\textwidth}
      \includegraphics[width=\linewidth]{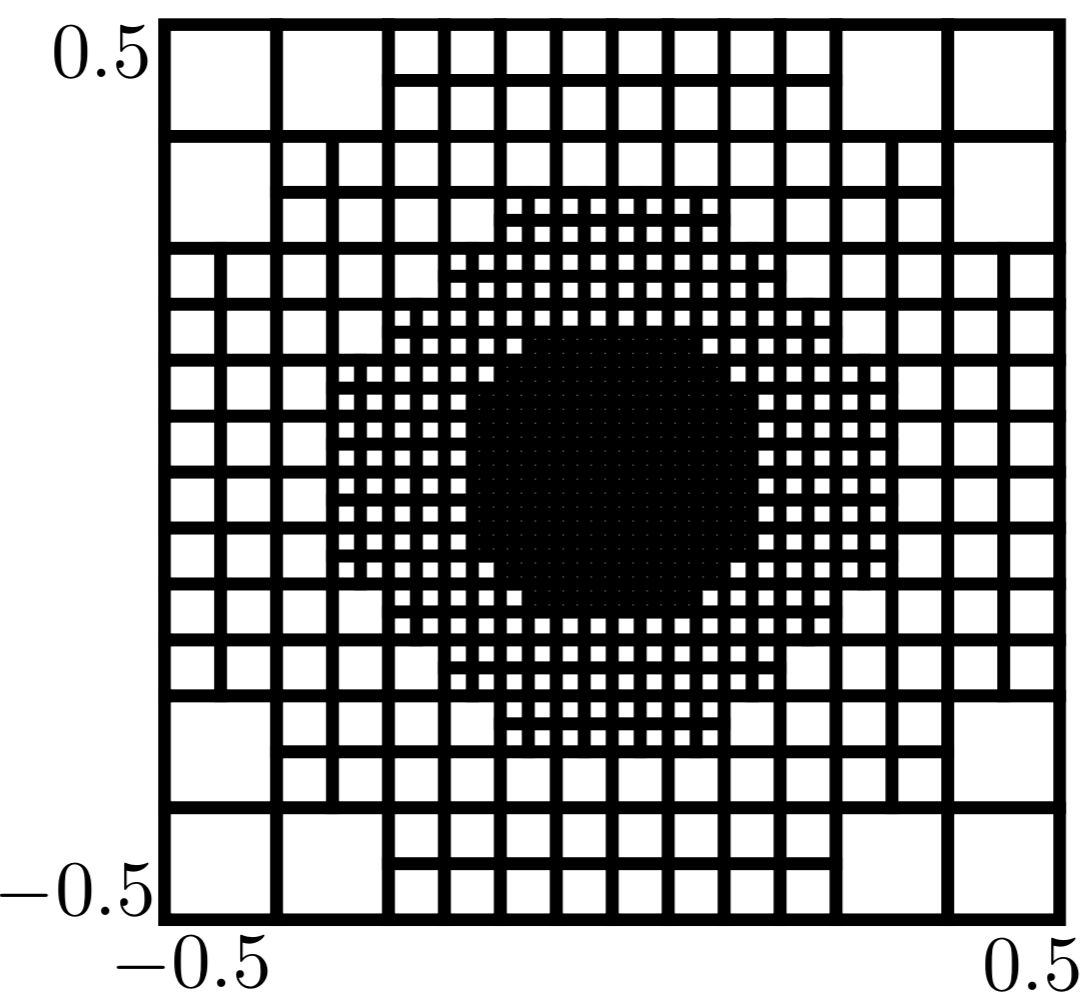}
        \caption{}
        \label{fig:grid_10}
    \end{subfigure}%
    \hfill
    \begin{subfigure}[b]{0.5\textwidth}
        \includegraphics[width=\linewidth]{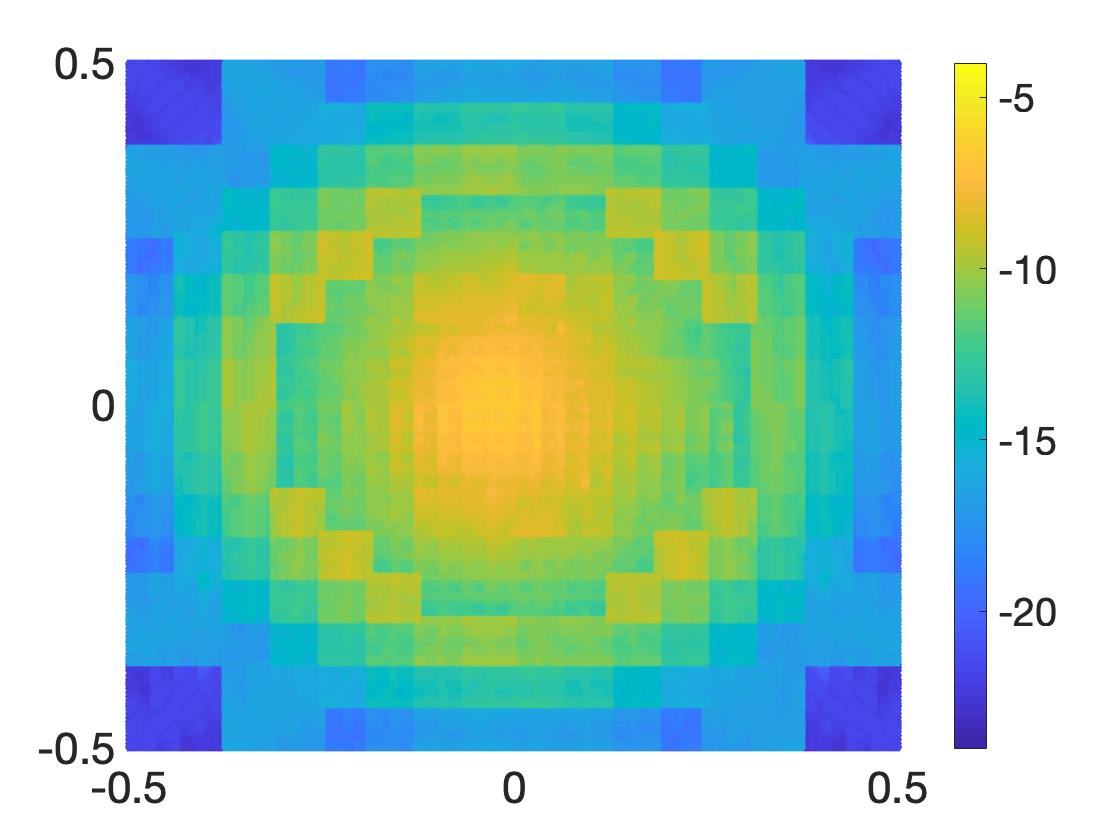}
        \caption{}
        \label{fig:error_grid_10}
    \end{subfigure}
    \caption{Results obtained with experiment 2. (a) and (b): The adaptive grid and log plot of error function with $\epsilon_{grid}=10^{-8}$ respectively; (c) and (d): The adaptive grid and log plot of error function with $\epsilon_{grid}=10^{-10}$ respectively.}
    \label{DAFMM_GMRES_40}
\end{center}
\end{figure}

We consider Gaussian contrast defined as
\begin{equation}\label{eq:gaussianContrast}
  q(x) = 1.5\exp(-160\left(x_1^2+x_2^2\right)).
\end{equation}
The leaf size, i.e., the number of grid points in a leaf, is set to $64$.
We set $\kappa$ to $40$, $\epsilon_{NCA}$ to $10^{-10}$ and $\epsilon_{GMRES}$ to $10^{-12}$. We solve for the scattered field, on a square $\Omega$, $[-0.5,0.5]^{2}$, using the DAFMM accelerated GMRES solver. To illustrate the convergence of the solver, we solve two systems generated with $\epsilon_{grid}$ set to $10^{-8}$ and $10^{-10}$.
With $\epsilon_{grid}=10^{-8}$, the generated system is of size $N=14848$. With $\epsilon_{grid}=10^{-10}$, $N=45568$.
The grids and log plot of error functions are given in Figure~\ref{DAFMM_GMRES_40}.
It is to be observed that the maximum value of the error function decreases as $\epsilon_{grid}$ decreases.
Plot of the Gaussian contrast and the real part of the field $u(x)$, obtained with $\epsilon_{grid}=10^{-10}$ are given in Figures~\ref{fig:GaussianContrast} and~\ref{fig:GMRES40Field} respectively.
The decay of residual with iteration count is shown in Figure~\ref{IterationTime40}. The CPU time to solve is shown in Figure~\ref{IterationTime40} and Table~\ref{table:IterationTime40}.
The plot of precomputation and apply time of DAFMM (to a single vector) versus $N$ is shown in Figure~\ref{DAFMM_single}, wherein we vary $N$ by varying $\epsilon_{grid}$.

\subsubsection{Experiment 3: DAFMM accelerated GMRES solver with HODLR preconditioner for Gaussian contrast with $\kappa=40$}
\label{Experiment2}
\begin{figure}[H]
  \begin{center}
    \begin{subfigure}[b]{0.5\textwidth}
          \includegraphics[width=\linewidth]{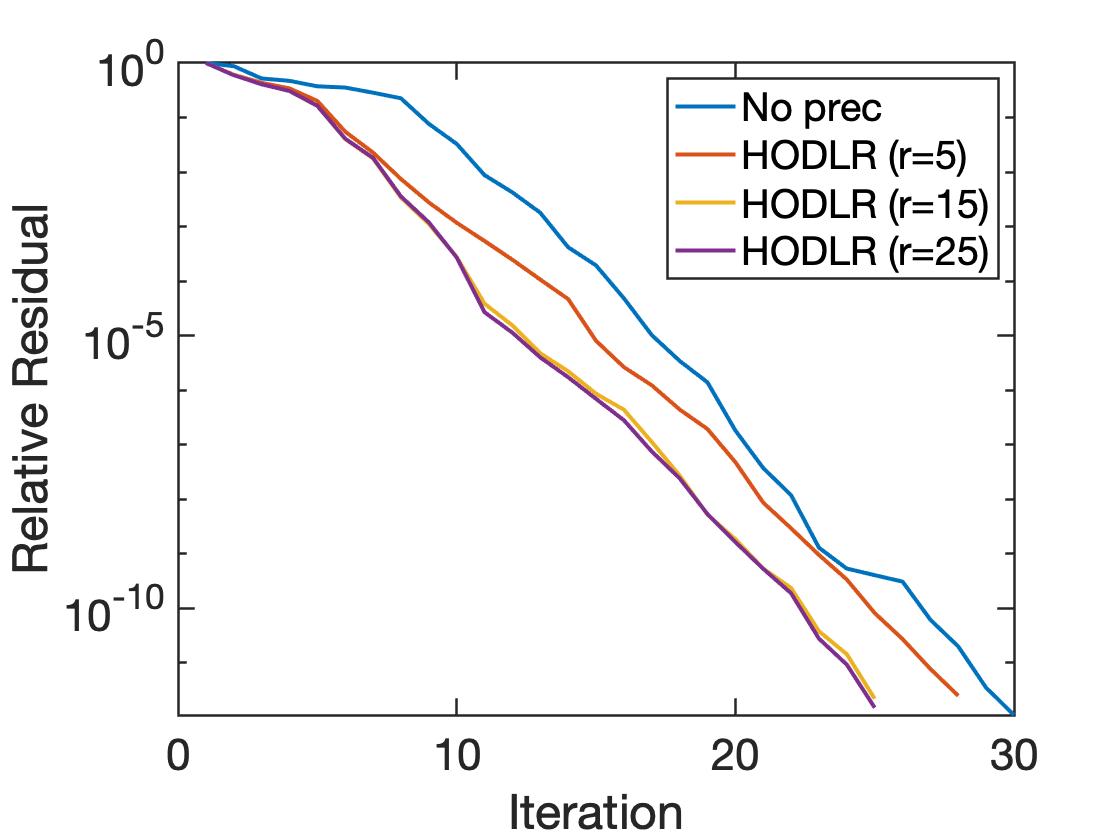}
          \label{fig:Iteration40}
      \end{subfigure}%
      \begin{subfigure}[b]{0.5\textwidth}
        \includegraphics[width=\linewidth]{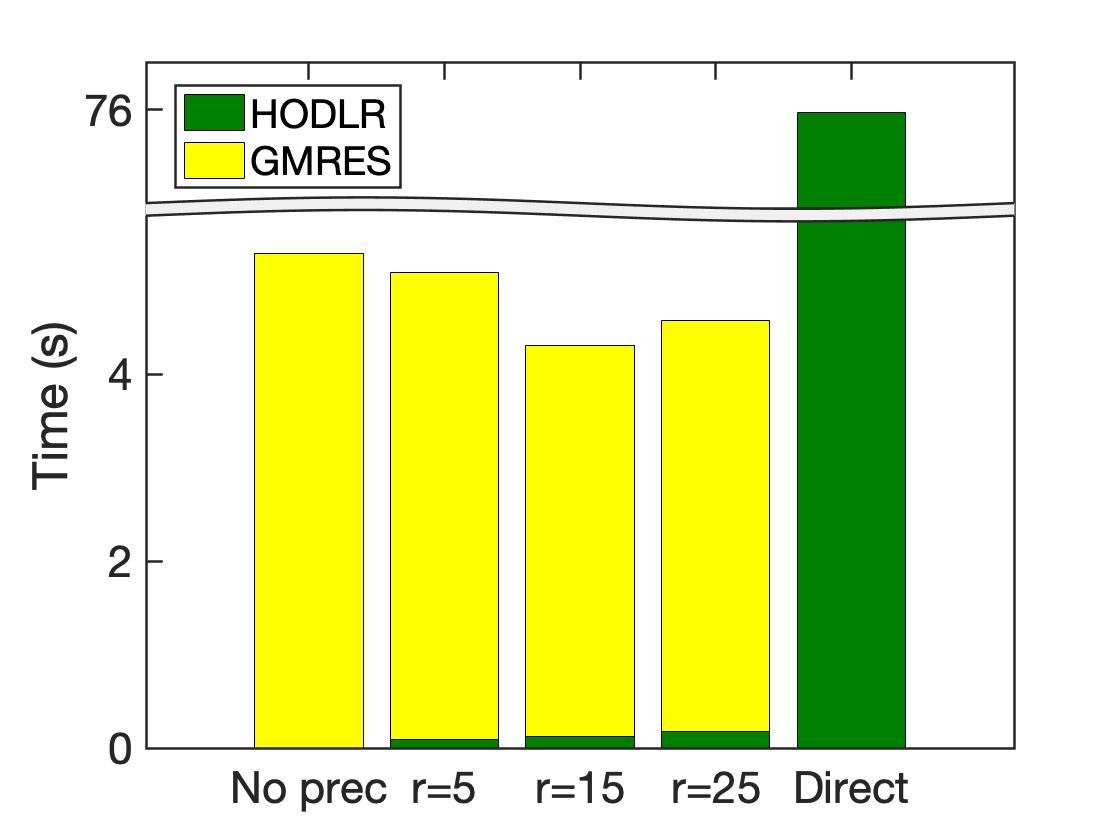}
          \label{fig:Time40}
      \end{subfigure}
      \caption{Results obtained with experiments 2,3 and 4 with $\epsilon_{grid}=10^{-8}$; `No prec' represents the GMRES solver with no pre-conditioner; `$r=5$', `$r=15$' and `$r=25$' represent Hybrid solvers where $r$ indicates the off-diagonal block rank of the HODLR pre-conditioner; `Direct' represents the HODLR direct solver; \textit{Left:} Relative residual $||A\hat{\psi}-f||_{2}/||f||_{2}$ as a function of the iteration count; \textit{Right:} Solve time. For the bars corresponding to Hybrid solvers, the green blocks indicate the time needed to factorise the preconditoner and the yellow blocks indicate the time needed to apply the pre-conditioner and solve by the GMRES solver.}
      \label{IterationTime40}
  \end{center}
\end{figure}

We set the leaf size to $64$, $\kappa$ to $40$, $\epsilon_{NCA}$ to $10^{-10}$, $\epsilon_{GMRES}$ to $10^{-12}$ and $\epsilon_{grid}$ to $10^{-8}$. We use Gaussian contrast as defined in equation \ref{eq:gaussianContrast}. The generated system is of size $N=14848$. We solve for the scattered field using the GMRES solver accelerated by DAFMM with HODLR as preconditioner on a square $\Omega$, $[-0.5,0.5]^{2}$.
The HODLR preconditioner is built such that its off-diagonal blocks are compressed to a user-specified rank. An increase in the rank leads to a more accurate preconditioner but would incur additional CPU time needed to build it. We studied this behavior by choosing ranks $5$, $15$, and $25$.
The decay of residual with iteration count is shown in Figure~\ref{IterationTime40}. The CPU times to factorize the HODLR preconditioner and solve (the time taken by the GMRES routine which includes the times to apply the preconditioner and DAFMM in each iteration) are shown in Figure~\ref{IterationTime40} and Table~\ref{table:IterationTime40}.
For the example considered, the choice of $15$ for the off-diagonal block rank gives the smallest computational time.

\begin{table}[H]
  \centering
  \setlength\tabcolsep{4pt}
    \begin{tabular}{|c|c|c|c|c|c|}
    \hline
     & No prec & $r=5$ & $r=15$ & $r=25$ & Direct \\ \hline\hline
     Factorization time (s) & $-$ & 0.087584 & 0.1269 & 0.171822 & 75.51\\ \hline
     Solve time (s) & 5.29146 & 4.99747 & 4.17626 & 4.4053 & 0.45\\ \hline
    \end{tabular}
    \caption{Results obtained with experiments 2,3 and 4 with $\epsilon_{grid}=10^{-8}$; The factorization time tabulated here for the Hybrid solvers indicates the time to factorize the preconditioner and that for the direct solver indicates the time to factorize the HODLR direct solver.}
    \label{table:IterationTime40}
 \end{table}

\subsubsection{Experiment 4: HODLR direct solver for Gaussian contrast with $\kappa=40$}
\begin{figure}[H]
  \begin{center}
    \begin{subfigure}[b]{0.4\textwidth}
        \includegraphics[width=\linewidth]{FiguresFinal/tree_grid_8_new.png}
        \label{fig:grid_8_hodlr}
    \end{subfigure}%
    \hfill
    \begin{subfigure}[b]{0.5\textwidth}
          \includegraphics[width=\linewidth]{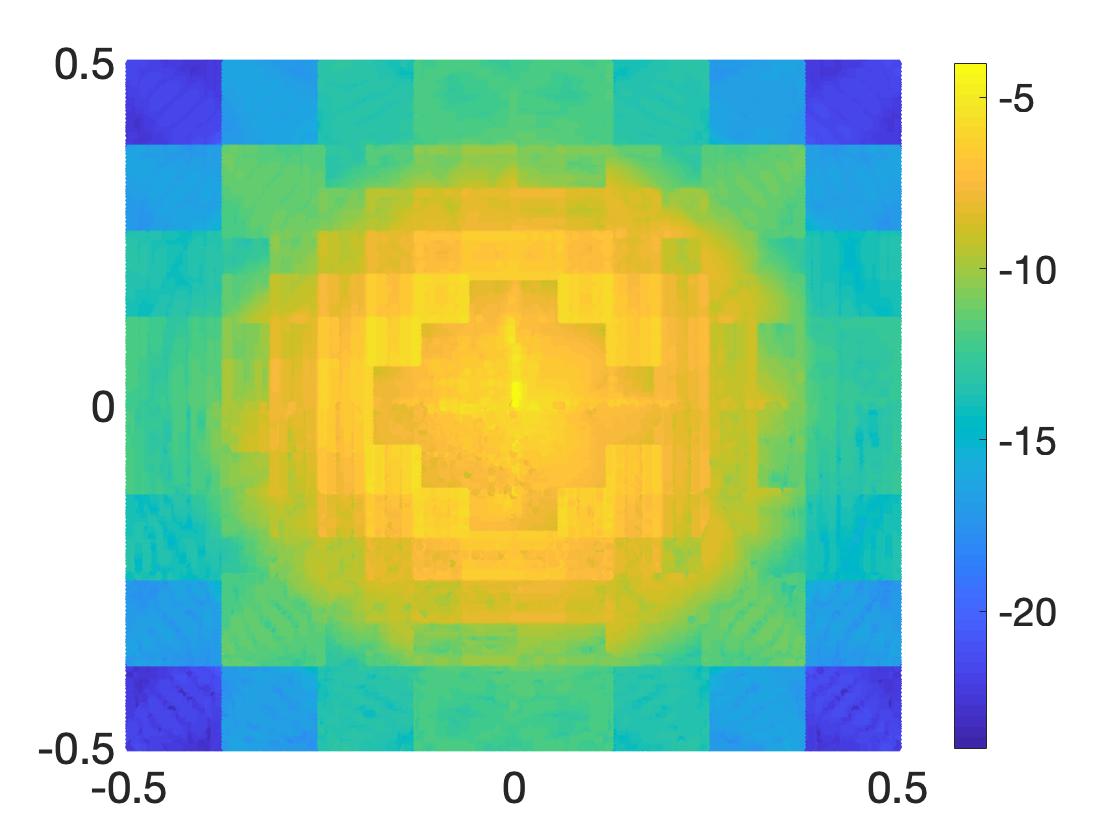}
          \label{fig:Iteration40}
      \end{subfigure}%
      \caption{Results obtained with experiment 4;  The adaptive grid and log plot of error function obtained with HODLR direct solver.}
      \label{HODLR_error}
  \end{center}
\end{figure}

We set the leaf size to $64$, $\kappa$ to $40$ and $\epsilon_{grid}$ to $10^{-8}$. We use Gaussian contrast as defined in equation~\eqref{eq:gaussianContrast}. The generated system is of size $N=14848$. We solve for the scattered field on a square $\Omega$, $[-0.5,0.5]^{2}$, using the HODLR direct solver, which is assembled such that the compression accuracy of the off-diagonal blocks is $10^{-10}$.
The grid and log plot of the error function are given in Figure~\ref{HODLR_error}.
The CPU times to factorize and solve are shown in Table~\ref{table:IterationTime40}. The sum of the factorization and solve times is plotted in Figure~\ref{IterationTime40} in comparison to the time taken by iterative solvers.
It is to be observed that the Hybrid solver with $r=15$ is more than $18$ times faster than the direct solver. This highlights the importance of the iterative solver and the HODLR pre-conditioner.
However, for the example considered, if one were to solve for roughly $20$ or more right hand sides, then HODLR direct solver is advantageous over the iterative solvers.

\subsubsection{Experiment 5: DAFMM accelerated GMRES solver with HODLR preconditioner for Gaussian contrast with $\kappa=300$} \label{Experiment3}
\begin{figure}[H]
  \begin{center}
    \begin{subfigure}[b]{0.5\textwidth}
            \includegraphics[width=\linewidth]{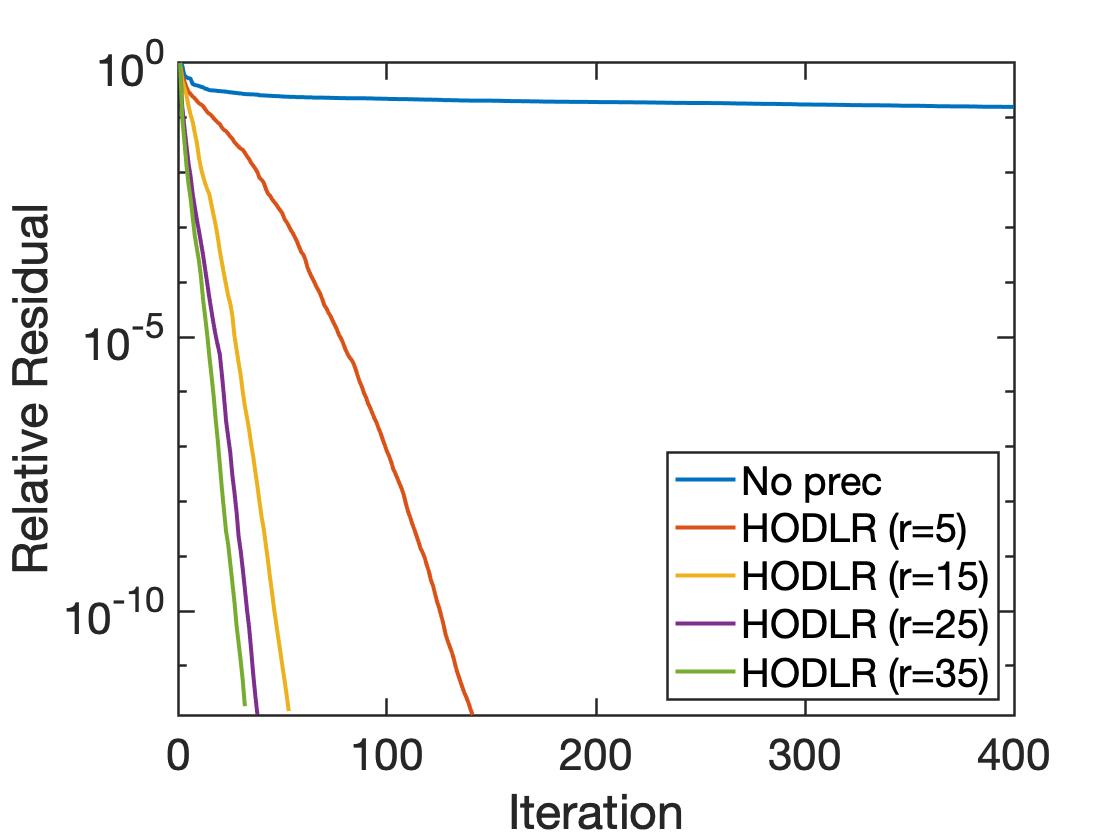}
            \label{fig:Iteration300}
        \end{subfigure}%
        \begin{subfigure}[b]{0.5\textwidth}
            \includegraphics[width=\linewidth]{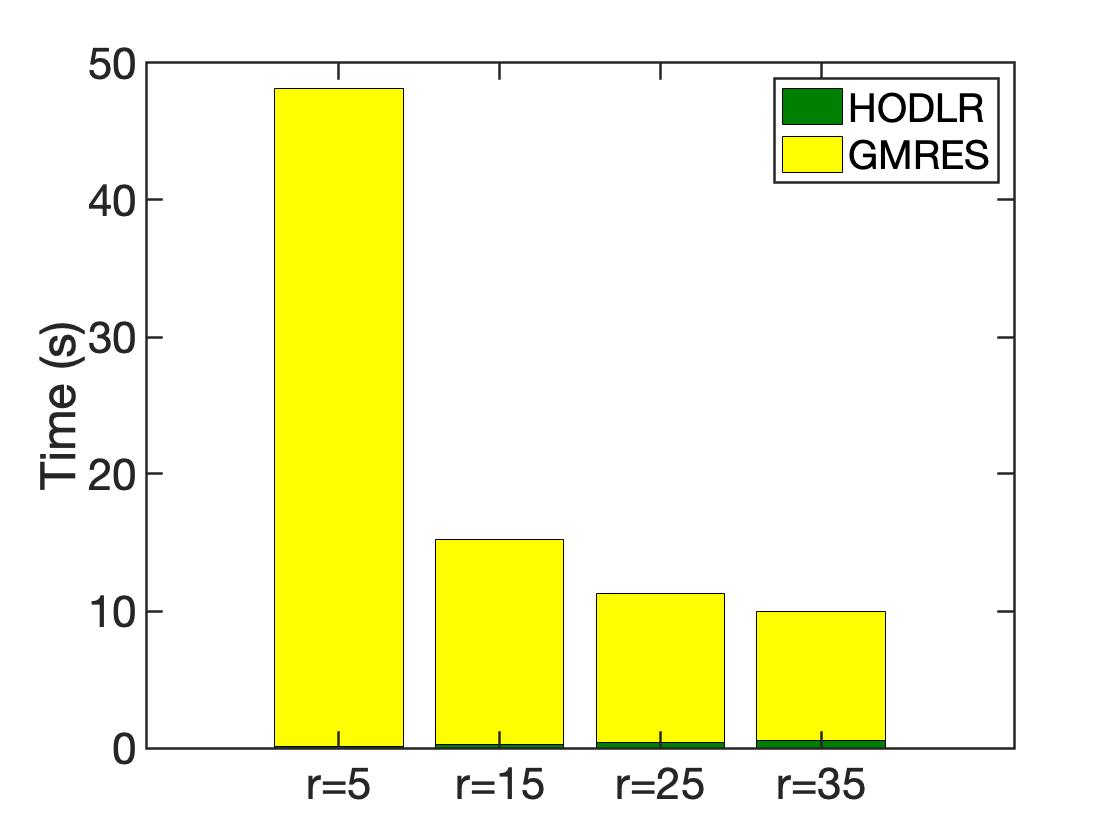}
            \label{fig:Time300}
        \end{subfigure}
        \caption{Results obtained with experiment 5; `No prec' represents the GMRES solver with no pre-conditioner; $r$ indicates the off-diagonal block rank of the HODLR pre-conditioner; \textit{Left:} Relative residual $||A\hat{\psi}-f||_{2}/||f||_{2}$ as a function of the iteration count; \textit{Right:} Solve time. The green blocks indicate the time needed to build the preconditoner. The yellow blocks indicate the time needed to apply the pre-conditioner and solve by the GMRES solver.}
        \label{IterationTime300}
\end{center}
\end{figure}

In this experiment we set $\kappa$ to $300$. We use Gaussian contrast as defined in equation~\eqref{eq:gaussianContrast}. $\epsilon_{grid}$ is set to $10^{-5}$.
The leaf size is set to $36$. With these inputs, the system generated is of size $36864$. To be conservative we set $\epsilon_{NCA}$ to $10^{-8}$ and $\epsilon_{GMRES}$ to $10^{-12}$.
We solve for the scattered field on a square $\Omega$, $[-0.5,0.5]^{2}$, using the GMRES solver accelerated by DAFMM with HODLR as a preconditioner.
 We experimented with different off-diagonal block ranks of the HODLR preconditioner.
The decay of residual with iteration count and the CPU time to solve are shown in Figure~\ref{IterationTime300}. For the example under consideration, the convergence of GMRES with no preconditioner was very slow. The relative residual after $400$ iterations is $0.158$. There is a significant improvement in convergence with the HODLR preconditioner, even with the rank of off-diagonal blocks set to $5$.
This highlights the importance of the HODLR preconditioner.

\begin{figure}[H]
  \begin{center}
    \begin{subfigure}[b]{0.4\textwidth}
            \includegraphics[width=\linewidth]{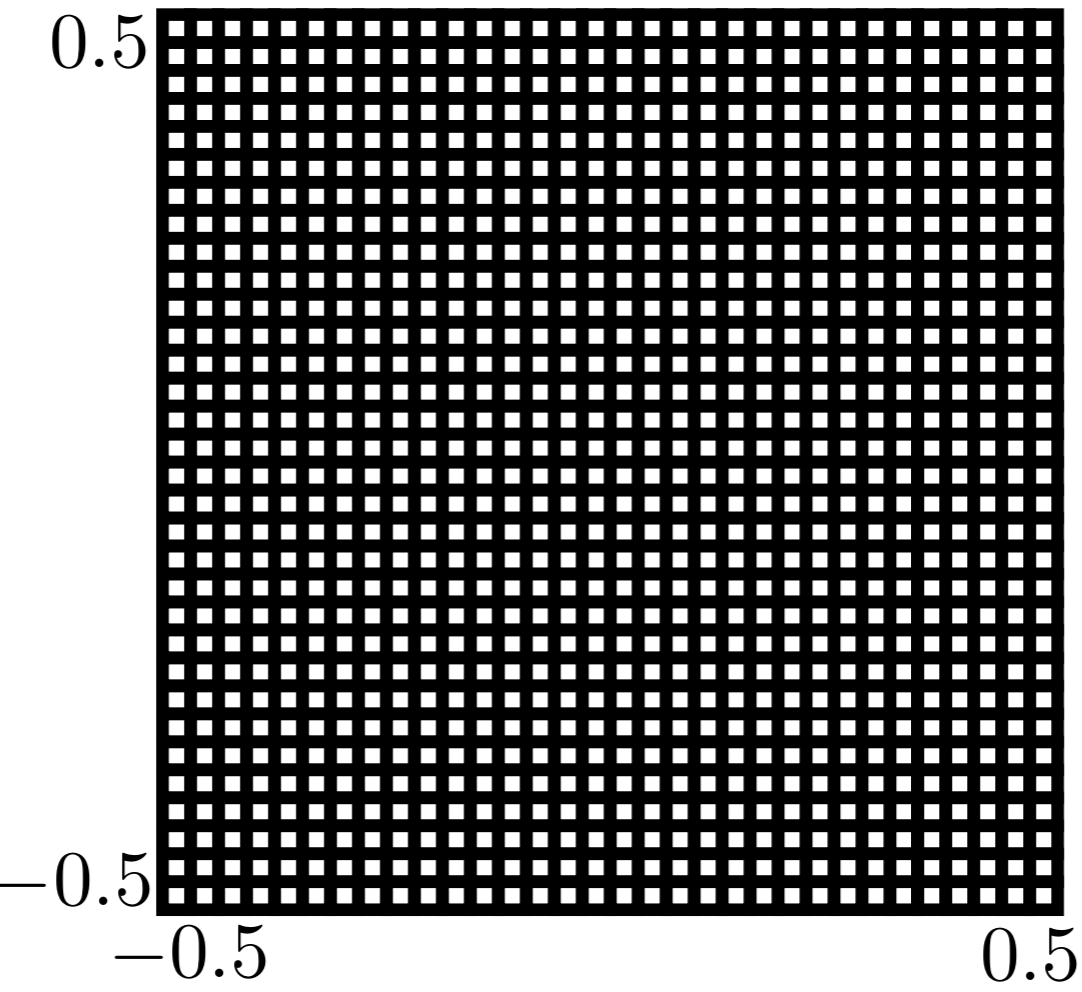}
            \caption{}
            \label{fig:tree300}
        \end{subfigure}%
        \begin{subfigure}[b]{0.5\textwidth}
            \includegraphics[width=\linewidth]{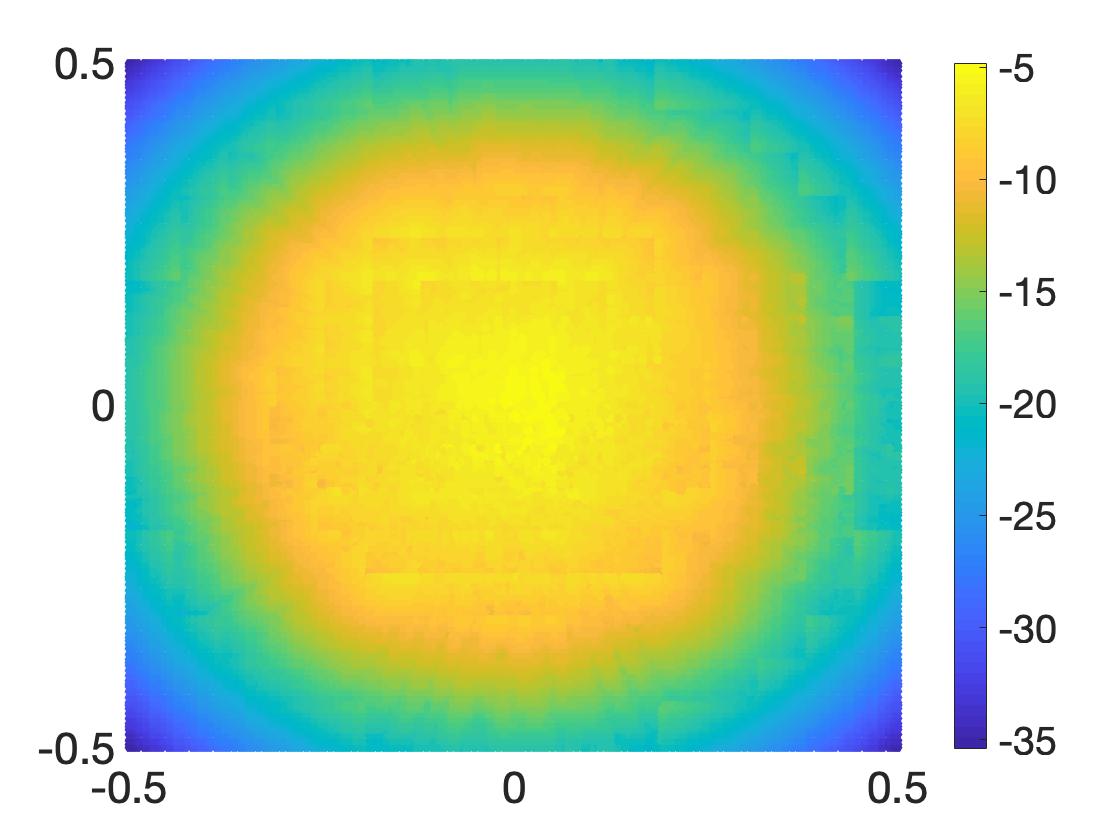}
            \caption{}
            \label{fig:ErrorPlot300}
        \end{subfigure}
        \caption{Results obtained with experiment 5: \textit{Left:} The adaptive grid; \textit{Right:} Log plot of error function}
        \label{DAFMM_GMRES_HODLR_300}
\end{center}
\end{figure}

With the off-diagonal block rank of HODLR set to $25$, we plot the real part of the field $u(x)$, the grid and the log plot of error function $E(x)$ in Figures~\ref{fig:GMRES300Field},~\ref{fig:tree300} and~\ref{fig:ErrorPlot300}.

\subsubsection{Experiment 6: Comparison of CPU times for direct and iterative solvers with Gaussian contrast}
We use Gaussian contrast as defined in equation~\eqref{eq:gaussianContrast}. We consider $\Omega$ to be the square $[-0.5,0.5]^{2}$. We set leaf size to $36$, $\epsilon_{GMRES}$ to $10^{-10}$ and $\kappa$ to $40$. The HODLR direct solver is assembled such that the compression accuracy of the off-diagonal blocks is $10^{-7}$. $\epsilon_{NCA}$ for the iterative solvers is also set to $10^{-7}$.
The HODLR preconditioner is assembled such that its off-diagonal blocks do not exceed a rank of $25$. We compare the CPU times of HODLR, GMRES, and Hybrid solvers in Table~\ref{table:gaussianTable} for various values of $N$, obtained by varying $\epsilon_{grid}$.
The decrease in speedup gained by Hybrid over GMRES solver with an increase in $N$, is due to the time incurred in constructing the pre-conditioner.

    \begin{table}[H]
      \centering
      \setlength\tabcolsep{4pt}
      \begin{tabular}{|c|c|c|c|c|c|c|c|c|c|c|}
        \hline
        & & \multicolumn{3}{|c|}{HODLR} & GMRES & \multicolumn{5}{|c|}{Hybrid}\\ \cline{3-11}
        N & $\epsilon_{grid}$ & T\textsubscript{Hf} & T\textsubscript{Hs} & T\textsubscript{HODLR} & T\textsubscript{GMRES} & T\textsubscript{Pf} & T\textsubscript{Ps} & T\textsubscript{Hybrid} & $\frac{\text{T\textsubscript{HODLR}}}{\text{T\textsubscript{Hybrid}}}$ & $\frac{\text{T\textsubscript{GMRES}}}{\text{T\textsubscript{Hybrid}}}$ \\ \hline\hline
        5328 & 1e-5 & 1.01 & 0.019 & 1.029 & 0.51 & 0.07 & 0.32 & 0.39 & 2.64  & 1.32 \\ \hline
        32112 & 1e-8 & 18.3 & 0.21 & 18.51 & 3.37 & 0.62 & 1.83 & 2.45 & 7.55 & 1.38 \\ \hline
        77904 & 1e-9 & 131.6 & 1.94 & 133.54 & 7.03 & 2.3 & 4.64 & 6.94 & 19.24 & 1.01 \\ \hline
        132768 & 1e-10 & 296.8 & 19.6 & 316.4 & 13.6 & 2.95 & 8.4 & 11.35 & 27.88 & 1.2 \\ \hline
        305568 & 1e-11 & 2548.9 & 176 & 2724.9 & 25.3 & 5.87 & 16.2 & 22.07 & 123.47 & 1.15 \\ \hline
        \end{tabular}
        \caption{Results obtained with experiment 6; CPU times of the three solvers for Gaussian contrast with $\kappa=40$.}
        \label{table:gaussianTable}
     \end{table}

\subsubsection{Experiment 7: Comparison of CPU times for direct and iterative solvers with multiple Gaussians contrast}
\begin{figure}[H]
  \begin{center}
    \begin{subfigure}[b]{0.33\textwidth}
        \includegraphics[width=\linewidth]{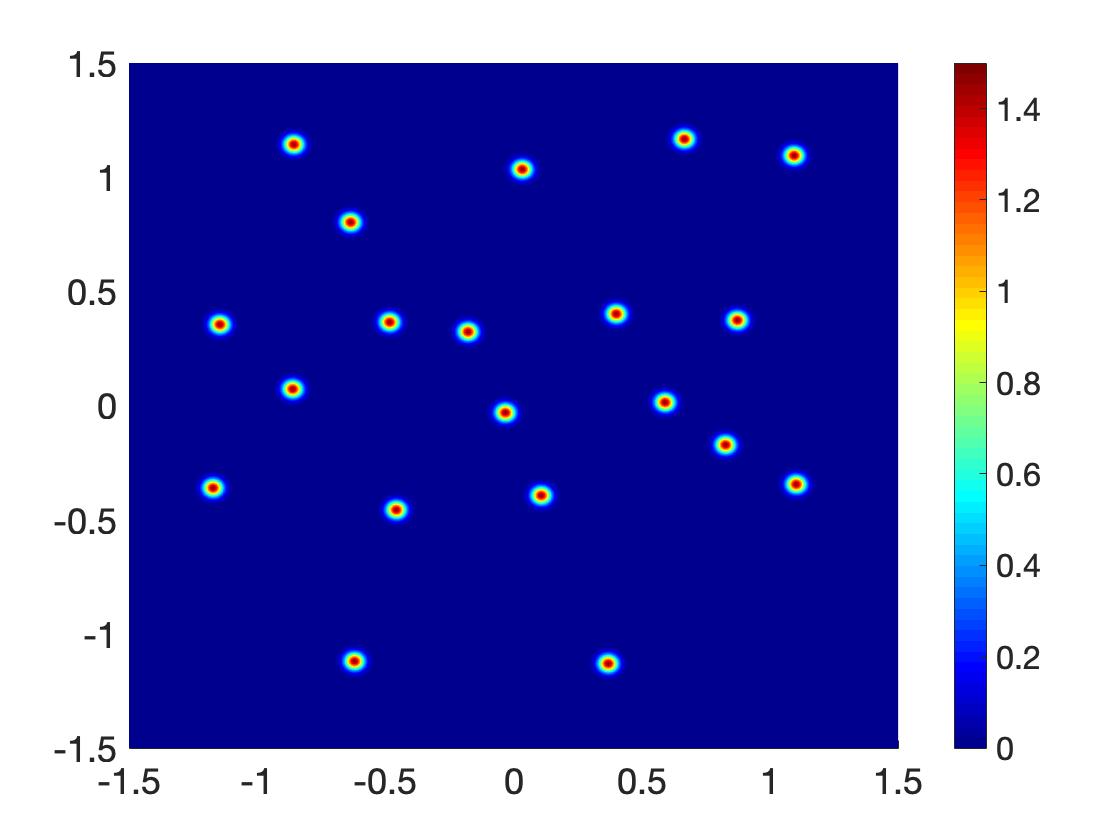}
        \caption{Contrast function}
        \label{fig:MultipleGaussianContrast}
    \end{subfigure}%
    \begin{subfigure}[b]{0.33\textwidth}
        \includegraphics[width=\linewidth]{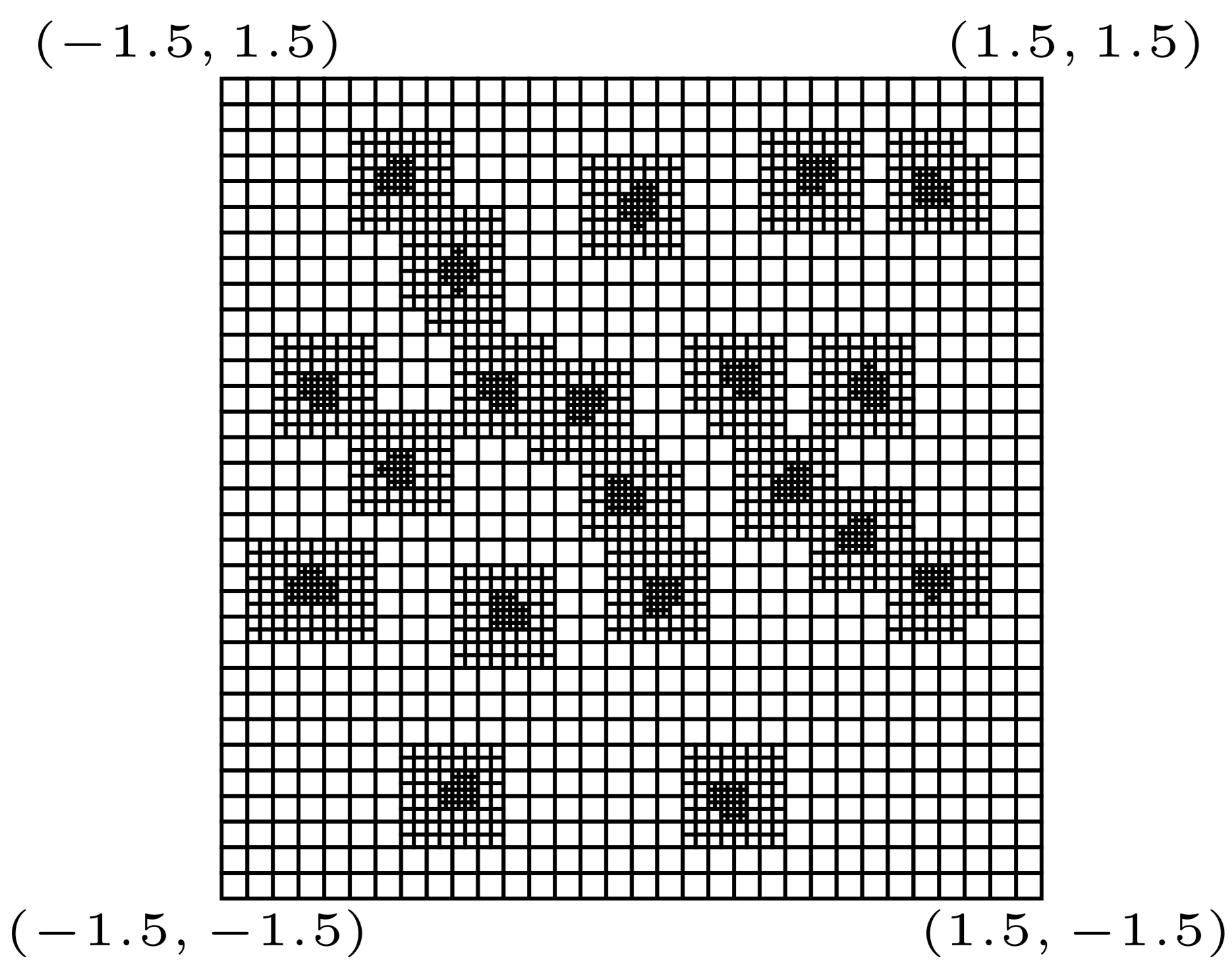}
        \caption{Adaptive Grid}
        \label{fig:MultipleGaussianGrid}
    \end{subfigure}%
    \begin{subfigure}[b]{0.33\textwidth}
        \includegraphics[width=\linewidth]{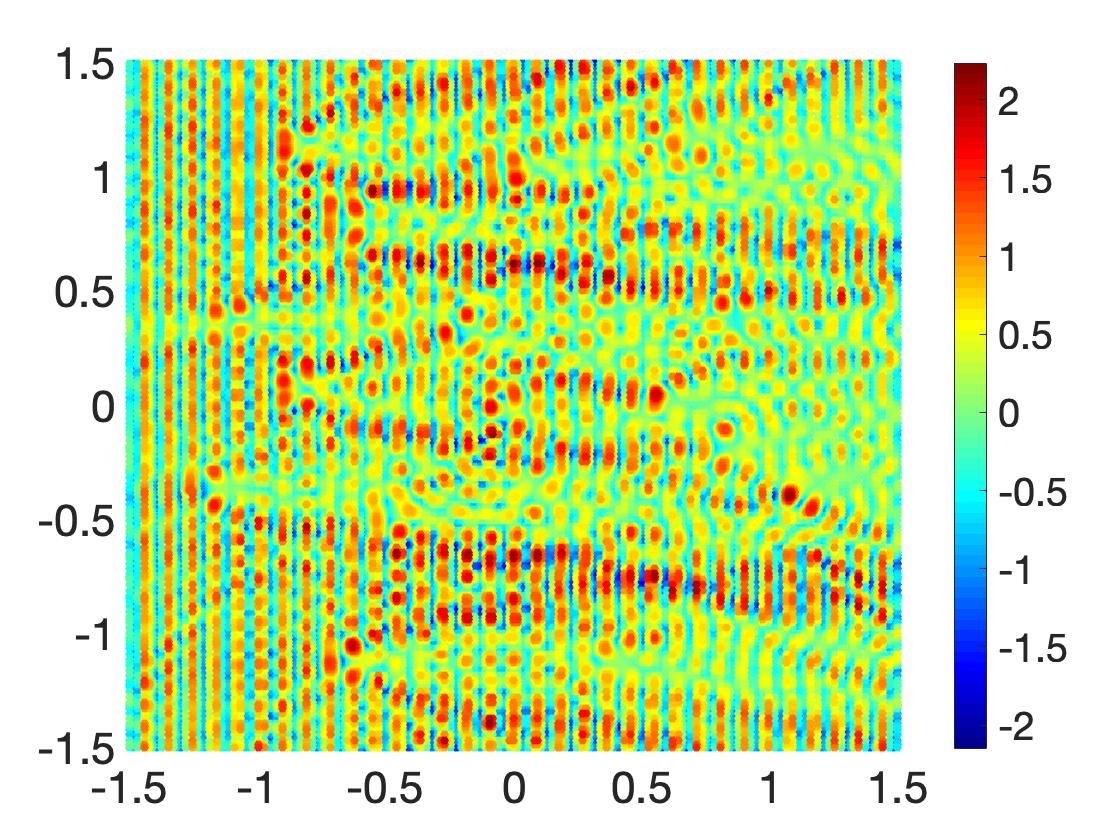}
        \caption{Real part of field}
        \label{fig:MultipleGaussianField}
    \end{subfigure}%
    \caption{Results obtained with experiment 7; Multiple Gaussian Contrast with $\kappa=70$ and $\epsilon_{grid}=10^{-4}$.}
    \label{fig:MultipleGaussianContrastPlots}
\end{center}
\end{figure}

We studied multiple scattering by simulating a scatterer with multiple Gaussians located in the domain. We choose $20$ well-separated Gaussians of the form
\begin{equation*}
  q_{i}(x) = 1.5\exp\left(\frac{-|x-x_{j}|^2}{a}\right)
\end{equation*}
centered randomly at $x_{j}$ in $[-1.5, 1.5]^{2}$, with $a=0.013$. Plot of the contrast function is illustrated in Figure~\ref{fig:MultipleGaussianContrast}.
We set leaf size to $36$, $\epsilon_{GMRES}$ to $10^{-10}$ and $\kappa$ to $70$. The HODLR direct solver is assembled such that the compression accuracy of the off-diagonal blocks is $10^{-7}$. $\epsilon_{NCA}$ for the iterative solvers is also set to $10^{-7}$.
The HODLR preconditioner is assembled such that its off-diagonal blocks do not exceed a rank of $15$. With $\epsilon_{grid}$ set to $10^{-4}$, we plot the adaptive grid and the real part of the total field evaluated using our iterative solver in Figures~\ref{fig:MultipleGaussianGrid} and~\ref{fig:MultipleGaussianField}. Table~\ref{table:multipleGaussiansTable} compares the solve time of the three solvers for different values of $\epsilon_{grid}$.

 \begin{table}[H]
   \centering
   \setlength\tabcolsep{4pt}
     \begin{tabular}{|c|c|c|c|c|c|c|c|c|c|c|}
     \hline
     & & \multicolumn{3}{|c|}{HODLR} & GMRES & \multicolumn{5}{|c|}{Hybrid}\\ \cline{3-11}
     N & $\epsilon_{grid}$ & T\textsubscript{Hf} & T\textsubscript{Hs} & T\textsubscript{HODLR} & T\textsubscript{GMRES} & T\textsubscript{Pf} & T\textsubscript{Ps} & T\textsubscript{Hybrid} & $\frac{\text{T\textsubscript{HODLR}}}{\text{T\textsubscript{Hybrid}}}$ & $\frac{\text{T\textsubscript{GMRES}}}{\text{T\textsubscript{Hybrid}}}$ \\ \hline\hline
     45828 & 1e-3 & 6.56 & 0.15 & 6.71 & 53.1 & 0.25 & 27.6 &  27.85 & 0.24 & 1.91 \\ \hline
    87948 & 1e-4 & 22.9 & 0.55 &  23.45 & 64.3 & 0.49 & 33.56 &  34.05 & 0.69 & 1.89 \\ \hline
    120240 & 1e-5 & 257.3 & 10.6 &  267.9 & 100.6 & 2.25 & 56.8 &  59.05 & 4.54  & 1.70 \\ \hline
    228564 & 1e-6 & 1572.5 & 101.1 &  1673.6 & 167.6 & 4.6 & 81.9 &  86.5 & 19.35 & 1.94 \\ \hline
    403524 & 1e-7 & 3116.2 & 347 &  3463.2 & 229.4 & 11.5 & 167.9 &  179.4 & 19.3 & 1.28 \\ \hline
    \end{tabular}
   \caption{Results obtained with experiment 7; CPU times of the three solvers for multiple Gaussians contrast.}
   \label{table:multipleGaussiansTable}
  \end{table}

\subsubsection{Experiment 8: Comparison of CPU times for direct and iterative solvers with Cavity contrast}
\begin{figure}[H]
  \begin{center}
    \begin{subfigure}[b]{0.33\textwidth}
        \includegraphics[width=\linewidth]{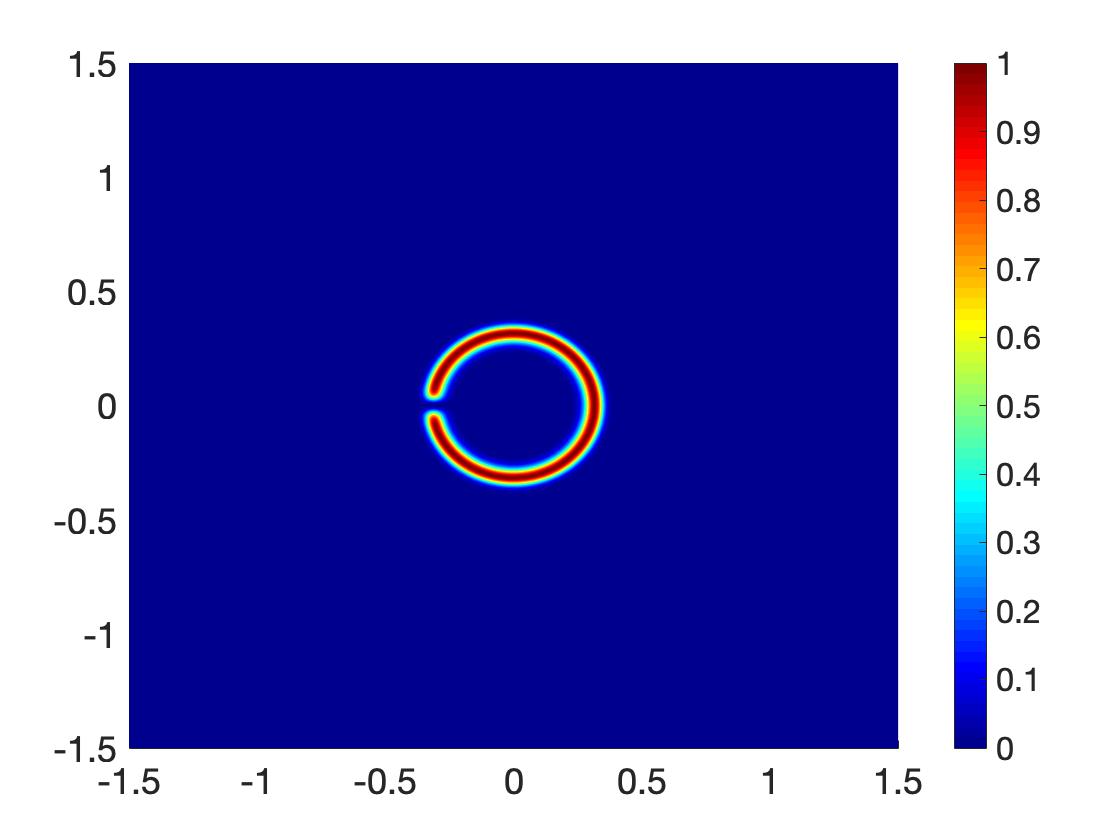}
        \caption{Contrast function}
        \label{fig:CavityContrast}
    \end{subfigure}%
    \begin{subfigure}[b]{0.33\textwidth}
        \includegraphics[width=\linewidth]{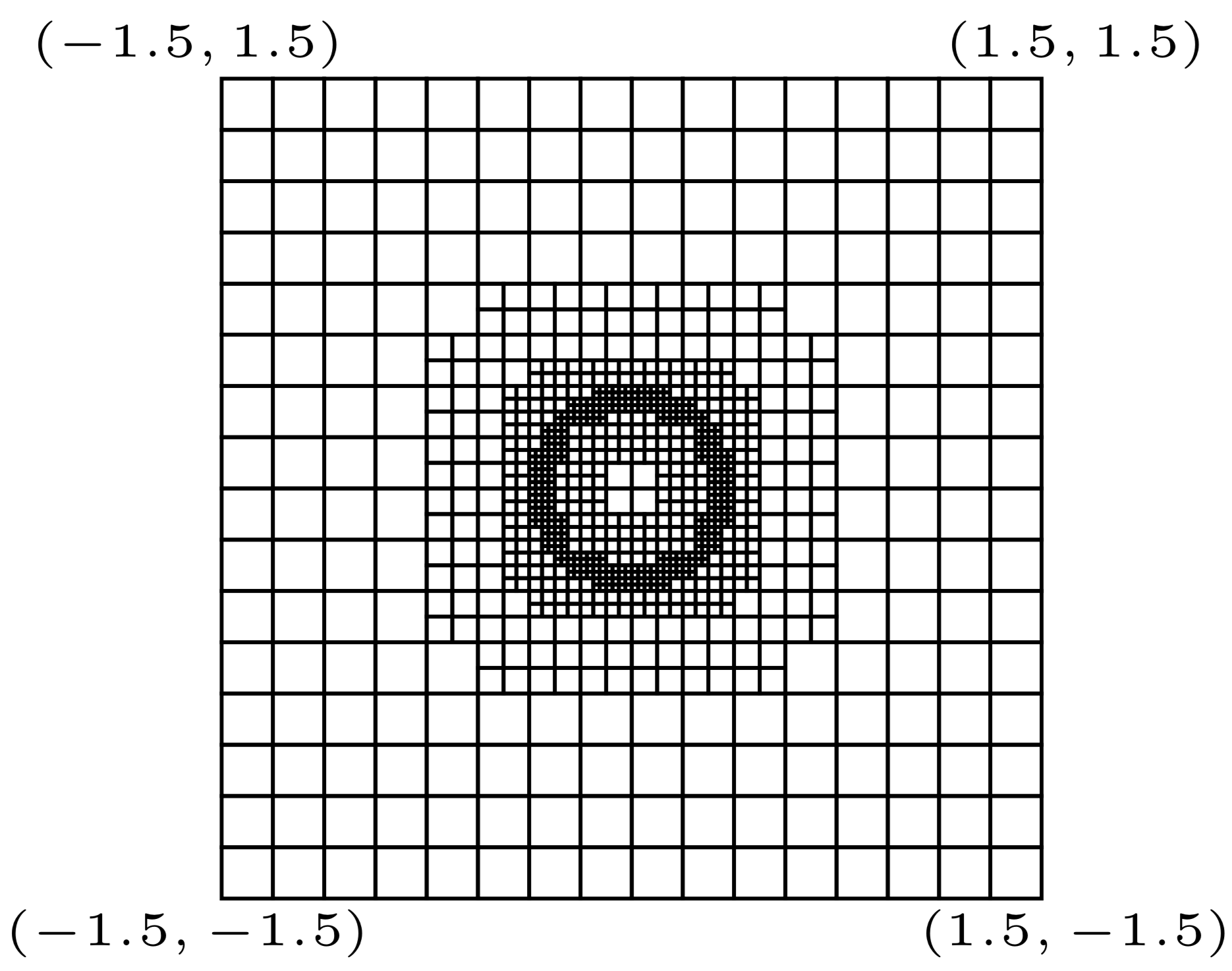}
        \caption{Adaptive Grid}
        \label{fig:CavityGrid}
    \end{subfigure}%
    \begin{subfigure}[b]{0.33\textwidth}
        \includegraphics[width=\linewidth]{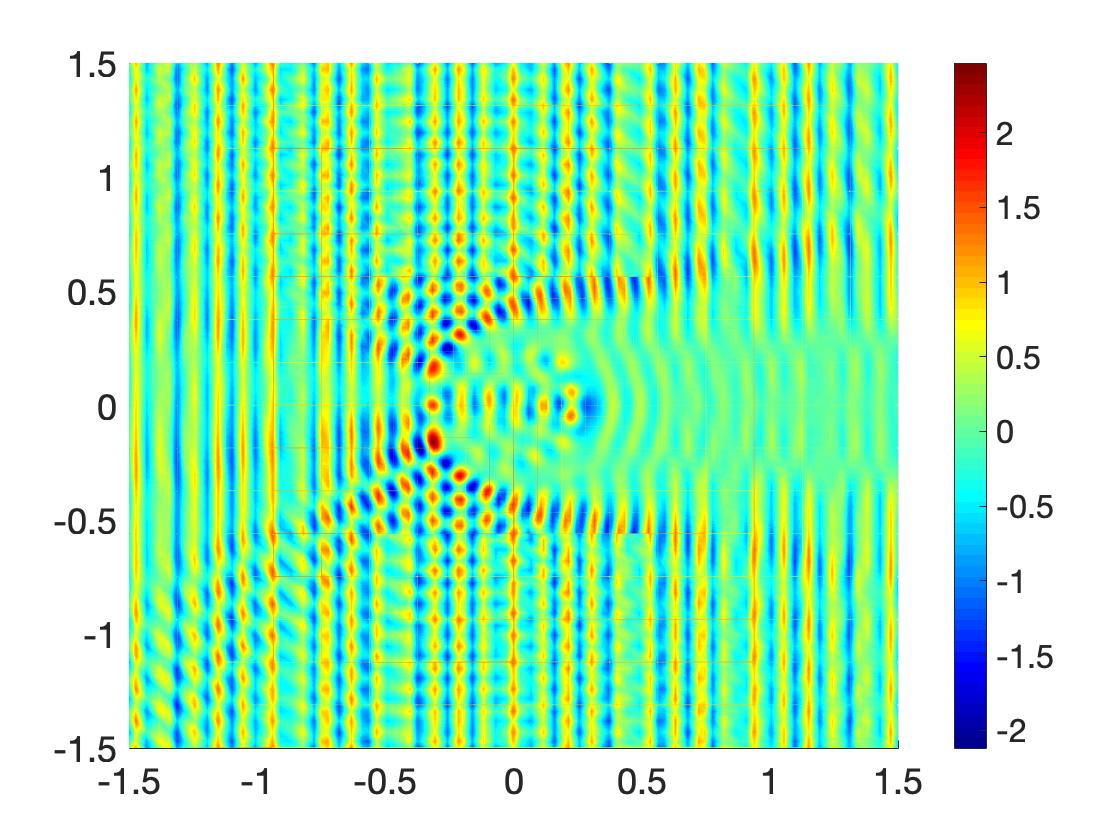}
        \caption{Real part of field}
        \label{fig:CavityField}
    \end{subfigure}%
    \caption{Results obtained with experiment 8; Cavity Contrast with $\kappa=60$ and $\epsilon_{grid}=10^{-4}$} 
    \label{fig:CavityContrastPlots}
\end{center}
\end{figure}

We consider the cavity contrast defined in polar coordinates as
$$q(r,\theta)= (1-\sin^{500}(0.5\theta))\exp(-2000(0.1-r^{2})^{2}).$$ Plot of the contrast function is illustrated in Figure~\ref{fig:CavityContrast}. We consider $\Omega$ to be the square $[-1.5,1.5]^{2}$.
We set leaf size to $36$, $\epsilon_{GMRES}$ to $10^{-10}$ and $\kappa$ to $60$. The HODLR direct solver is assembled such that the compression accuracy of the off-diagonal blocks is $10^{-7}$. $\epsilon_{NCA}$ for the iterative solvers is also set to $10^{-7}$.
The HODLR preconditioner is assembled such that its off-diagonal blocks do not exceed a rank of $25$. With $\epsilon_{grid}$ set to $10^{-4}$, we plot the adaptive grid and the real part of the total field evaluated using our iterative solver in Figures~\ref{fig:CavityGrid} and~\ref{fig:CavityField}. Table~\ref{table:cavityTable} compares the performances of the three solvers for different values of $\epsilon_{grid}$.

 \begin{table}[H]
   \centering
   \setlength\tabcolsep{4pt}
     \begin{tabular}{|c|c|c|c|c|c|c|c|c|c|c|}
     \hline
     & & \multicolumn{3}{|c|}{HODLR} & GMRES & \multicolumn{5}{|c|}{Hybrid}\\ \cline{3-11}
     N & $\epsilon_{grid}$ & T\textsubscript{Hf} & T\textsubscript{Hs} & T\textsubscript{HODLR} & T\textsubscript{GMRES} & T\textsubscript{Pf} & T\textsubscript{Ps} & T\textsubscript{Hybrid} & $\frac{\text{T\textsubscript{HODLR}}}{\text{T\textsubscript{Hybrid}}}$ & $\frac{\text{T\textsubscript{GMRES}}}{\text{T\textsubscript{Hybrid}}}$ \\ \hline\hline
     17856 & 1e-3 & 6.1 & 0.11 & 6.21 & 12.6 & 0.28 & 6.4 &  6.68 & 0.93 & 1.88 \\ \hline
     35568 & 1e-4 & 9.77 & 0.2 &  9.97 & 28.9 & 0.53 & 13.1 &  13.63 & 0.73 & 2.12 \\ \hline
     47232 & 1e-5 & 27.4 & 0.27 &  27.67 & 32.3 & 1.27 & 17.9 &  19.17 & 1.44 & 1.68 \\ \hline
     101880 & 1e-6 & 194.3 & 4.8 &  199.1 & 109.7 & 3.31 & 36 &  39.31 & 5.06 & 2.79 \\ \hline
     145080 & 1e-7 & 384.3 & 35.5 &  419.8 & 156.9 & 5.76 & 48.8 &  54.56 & 7.69 & 2.88 \\ \hline
     362376 & 1e-8 & 3889 & 451 &  4340 & 450 & 35.87 & 149 & 184.87 & 23.48 & 2.43 \\ \hline
     \end{tabular}
     \caption{Results obtained with experiment 8; CPU times of the three solvers for cavity contrast.}
     \label{table:cavityTable}
  \end{table}

\subsubsection{Experiment 9: Comparison of CPU times for direct and iterative solvers with Lens contrast}
\label{Experiment7}
\begin{figure}[H]
  \begin{center}
    \begin{subfigure}[b]{0.33\textwidth}
        \includegraphics[width=\linewidth]{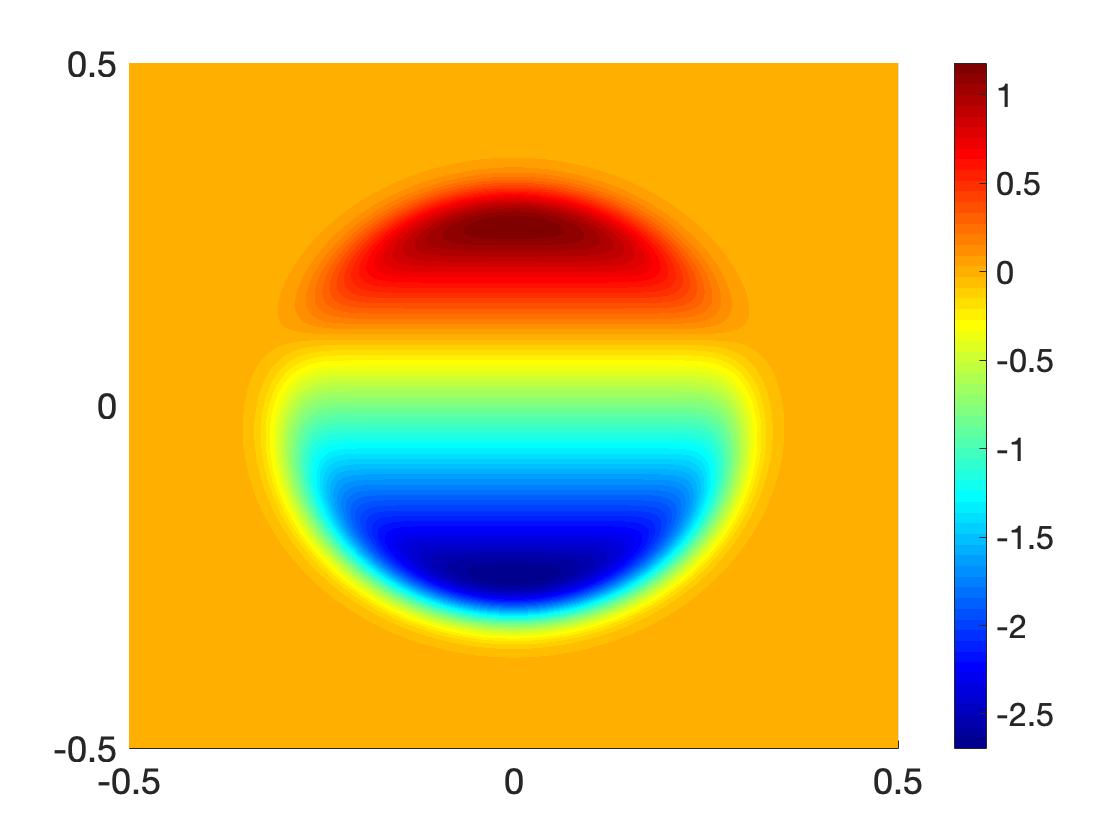}
        \caption{Contrast function}
        \label{fig:LensContrast}
    \end{subfigure}%
    \begin{subfigure}[b]{0.26\textwidth}
        \includegraphics[width=\linewidth]{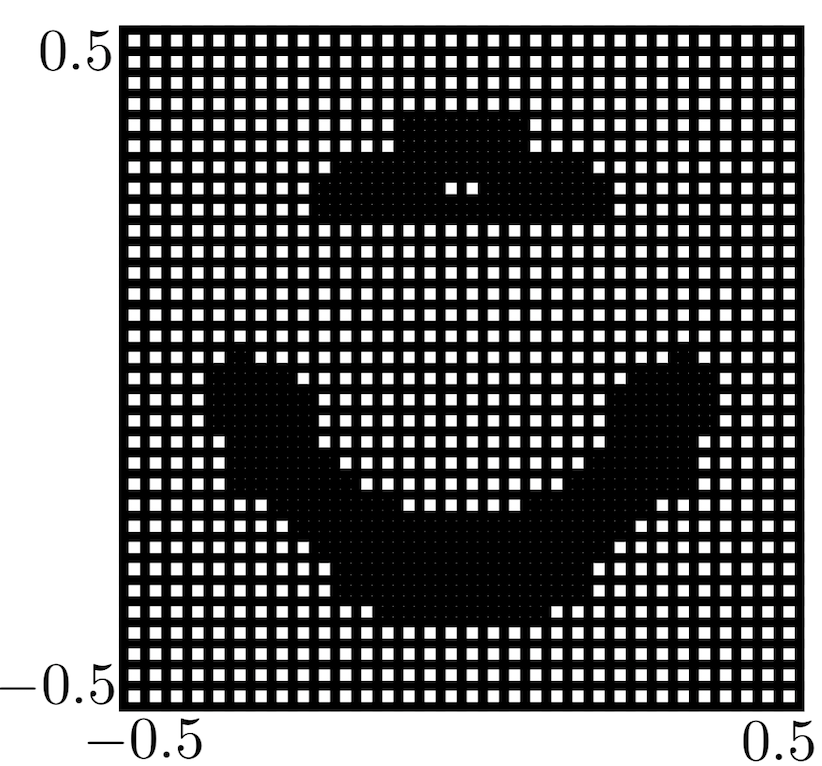}
        \caption{Adaptive Grid}
        \label{fig:LensGrid}
    \end{subfigure}%
    \begin{subfigure}[b]{0.33\textwidth}
        \includegraphics[width=\linewidth]{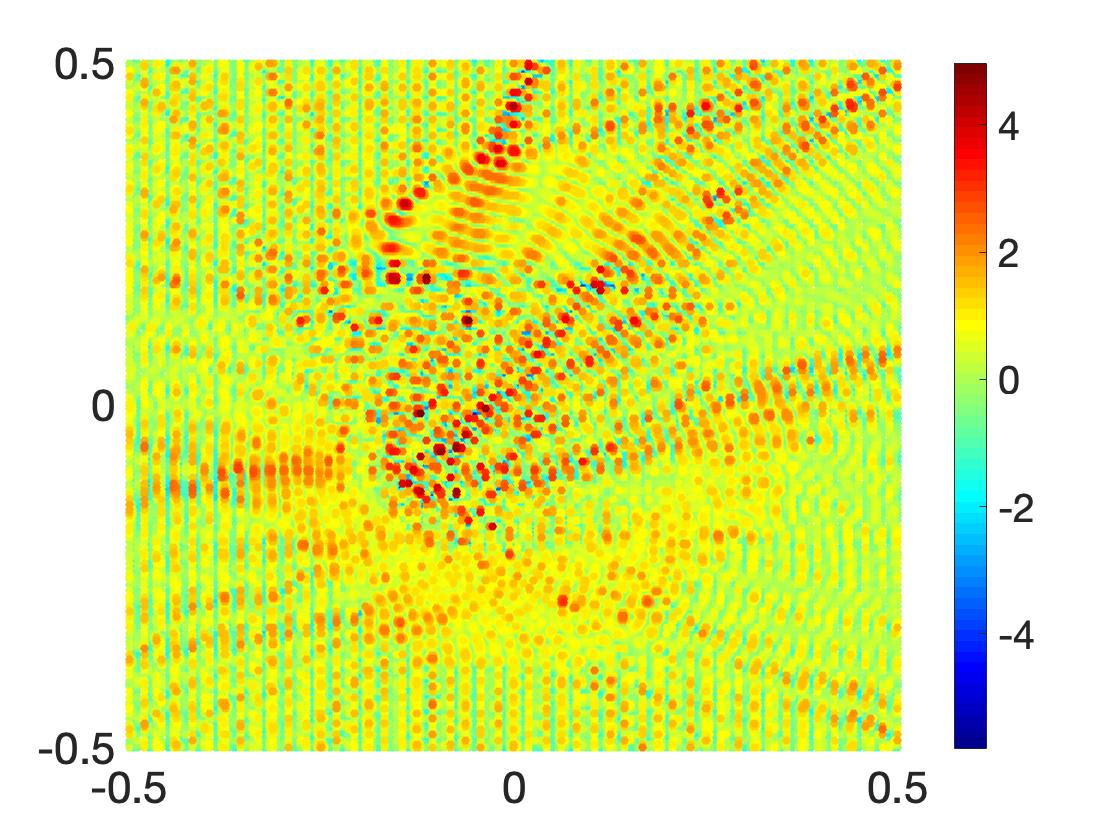}
        \caption{Real part of field}
        \label{fig:LensField}
    \end{subfigure}%
    \caption{Results obtained with experiment 9; Lens Contrast with $\kappa=300$ and $\epsilon_{grid}=10^{-6}$.}
    \label{fig:LensContrastPlots}
\end{center}
\end{figure}

We consider the lens contrast defined as
$$q(x)= 4(x_2-0.1)(1-\text{erf}(25((x_1^{2}+x_2^{2})^{0.5}-0.3))).$$
Plot of the contrast function is illustrated in Figure~\ref{fig:LensContrast}. We consider $\Omega$ to be the square $[-0.5,0.5]^{2}$.
We set leaf size to $36$, $\epsilon_{GMRES}$ to $10^{-10}$ and $\kappa$ to $300$. The HODLR direct solver is assembled such that the compression accuracy of the off-diagonal blocks is $10^{-7}$. $\epsilon_{NCA}$ for the iterative solvers is also set to $10^{-7}$.
The HODLR preconditioner is assembled such that its off-diagonal blocks do not exceed a rank of $70$. With $\epsilon_{grid}$ set to $10^{-6}$, we plot the adaptive grid and the real part of the total field evaluated using our iterative solver in Figures~\ref{fig:LensGrid} and~\ref{fig:LensField}. Table~\ref{table:lensTable} compares the performances of HODLR and Hybrid solvers for different values of $\epsilon_{grid}$.
The GMRES timing is not reported as the GMRES without preconditioning takes thousands of iterations to converge (even for \textit{smaller} values of $N$).

 \begin{table}[H]
   \centering
   \setlength\tabcolsep{4pt}
     \begin{tabular}{|c|c|c|c|c|c|c|c|c|}
     \hline
     & & \multicolumn{3}{|c|}{HODLR} & \multicolumn{4}{|c|}{Hybrid}\\ \cline{3-9}
     N & $\epsilon_{grid}$ & T\textsubscript{Hf} & T\textsubscript{Hs} & T\textsubscript{HODLR} & T\textsubscript{Pf} & T\textsubscript{Ps} & T\textsubscript{Hybrid} & $\frac{\text{T\textsubscript{HODLR}}}{\text{T\textsubscript{Hybrid}}}$ \\ \hline\hline
     36864 & 1e-5 & 46.2 & 0.39 & 46.59 & 3.2 & 29 &  32.2 & 1.45 \\ \hline
     56736 & 1e-6 & 196.6 & 1.18 &  197.78 & 12.25 & 98.9 &  111.1 & 1.78 \\ \hline
     83520 & 1e-7 & 374.1 & 27.16 &  401.26 & 27.25 & 182.1 &  209.3 & 1.92 \\ \hline
     182664 & 1e-8 & 2512.5 & 162.51 &  2675.01 & 100.1 & 1128.8 &  1228.9 & 2.18 \\ \hline
     \end{tabular}
 \caption{Results obtained with experiment 9; CPU times of HODLR and Hybrid solvers for lens contrast.}
     \label{table:lensTable}
  \end{table}

 \begin{figure}[H]
  \begin{subfigure}[t]{.49\textwidth}
    \centering
    \includegraphics[width=\linewidth]{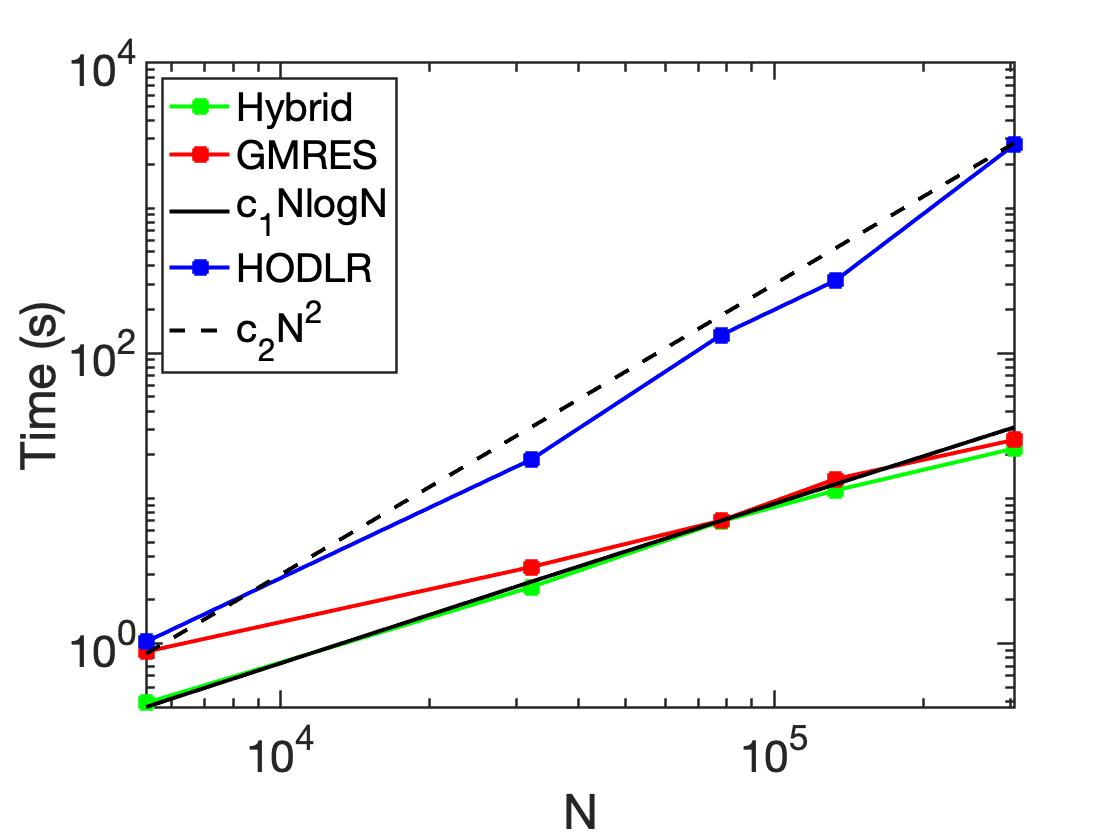}
    \caption{Gaussian Contrast}
    \label{fig:gaussianScaling}
  \end{subfigure}
  \hfill
  \begin{subfigure}[t]{.49\textwidth}
    \centering
    \includegraphics[width=\linewidth]{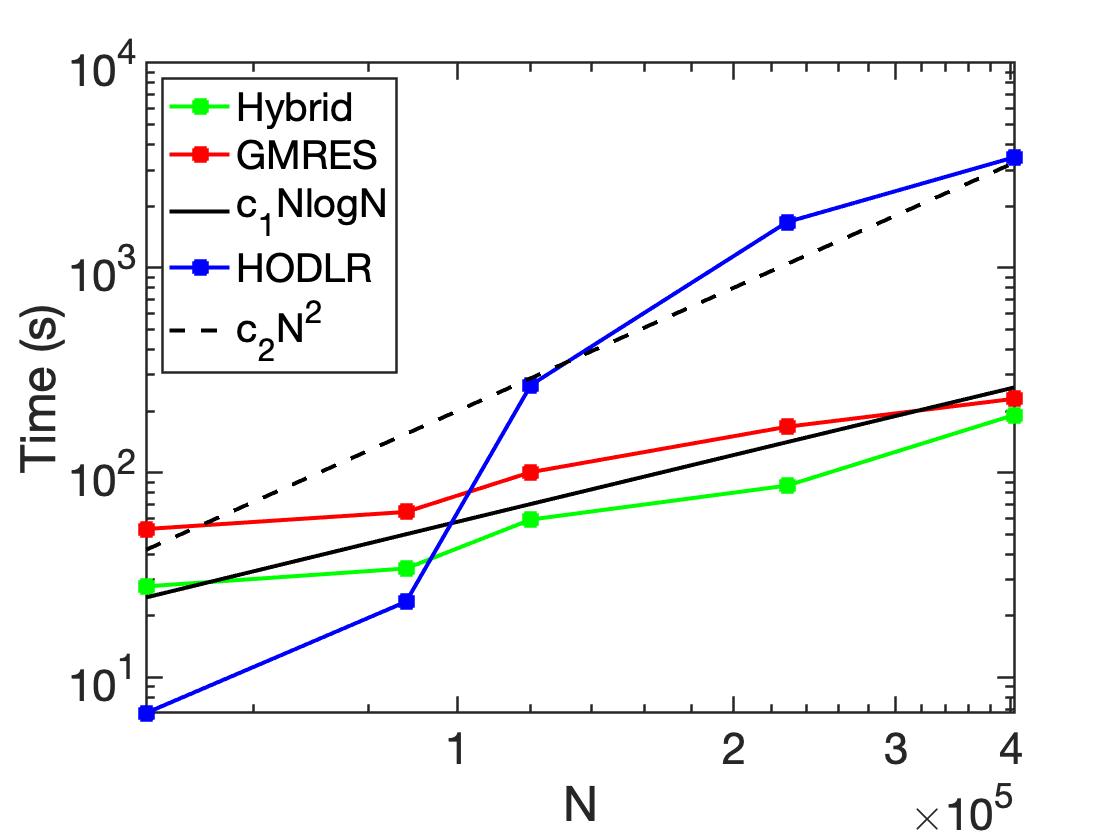}
    \caption{Multiple Gaussians Contrast}
    \label{fig:multipleGaussiansScaling}
  \end{subfigure}
  \medskip

  \begin{subfigure}[t]{.49\textwidth}
    \centering
    \includegraphics[width=\linewidth]{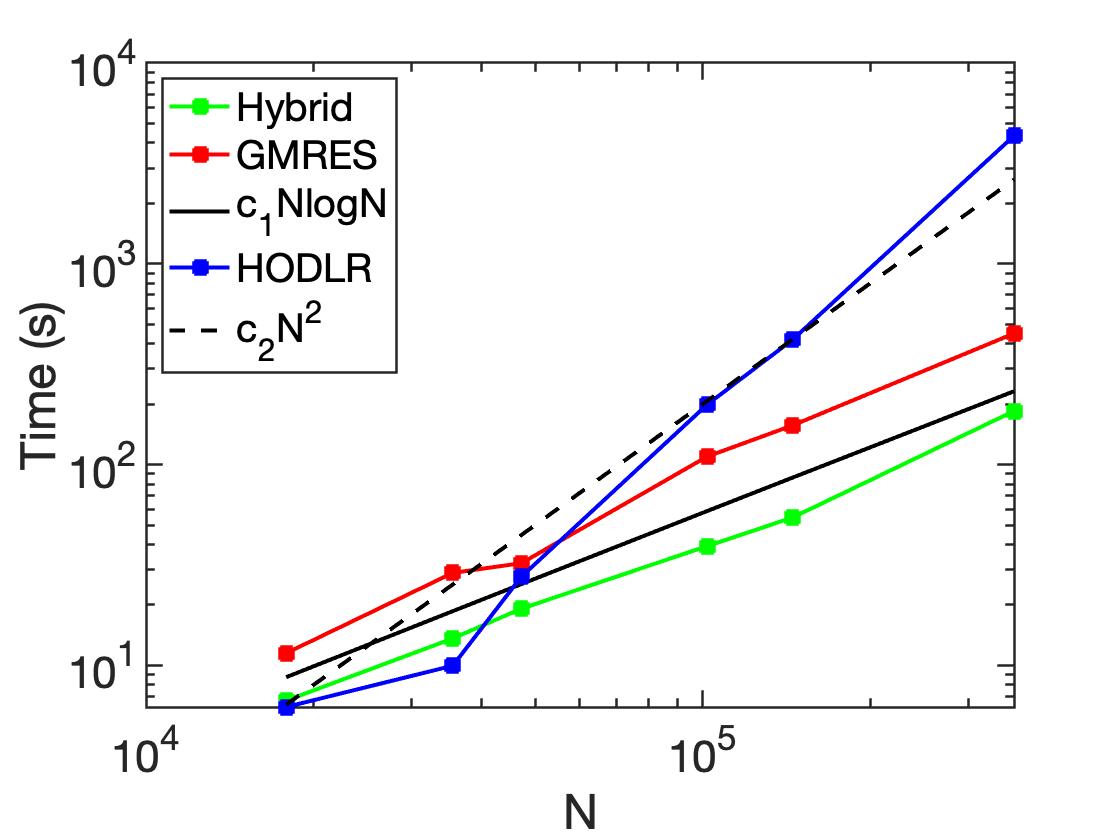}
    \caption{Cavity Contrast}
    \label{fig:cavityScaling}
  \end{subfigure}
  \hfill
  \begin{subfigure}[t]{.49\textwidth}
    \centering
    \includegraphics[width=\linewidth]{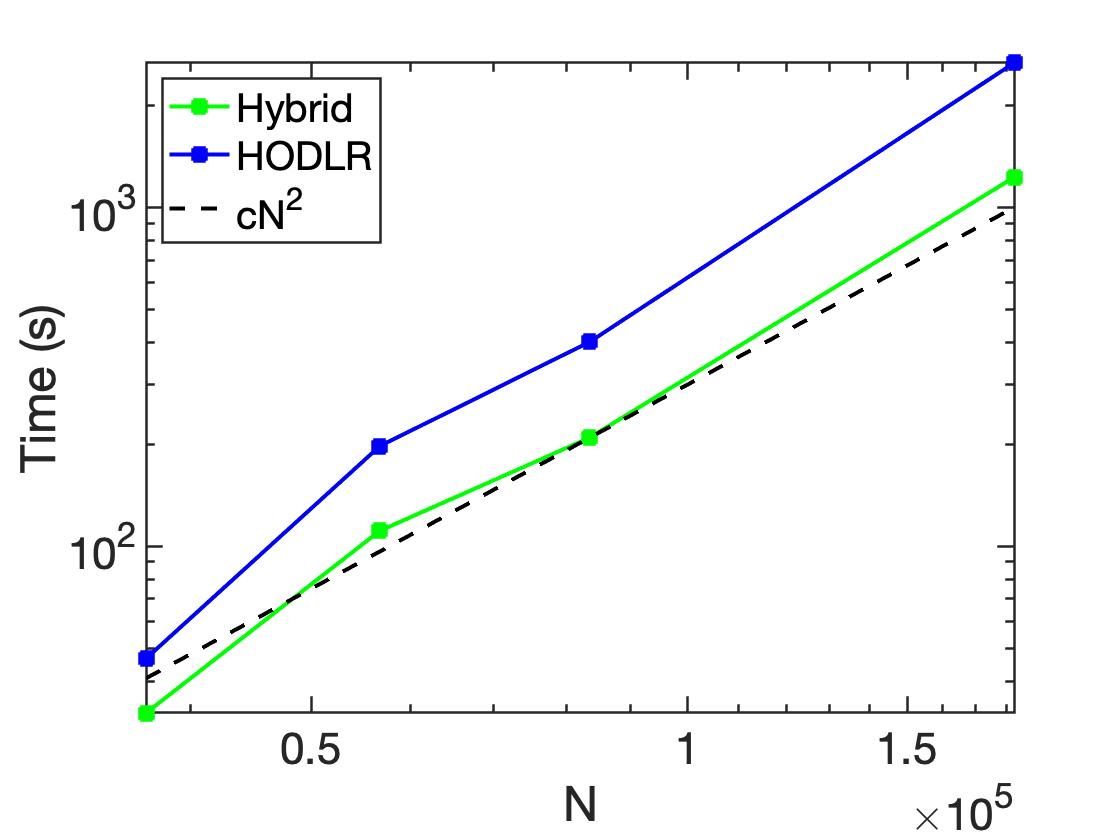}
    \caption{Lens Contrast}
    \label{fig:lensScaling}
  \end{subfigure}
  \caption{Results obtained with experiments 6-9; Time complexities of the different solvers for different contrast functions}
  \label{fig:Scaling}
\end{figure}

\subsection{Inferences}
Sections~\ref{Experiment3} and~\ref{Experiment7} indicate that at high values of $\kappa$ (the example considered is $300$), the convergence of GMRES is poor.
Hence in such cases, it is advisable to either use an iterative solver with a preconditioner or a direct solver. One possible reason could be the ill-conditioned nature of the linear system at high frequencies.

We observe that, in general, the hybrid solver outperforms the fast iterative solver accelerated by DAFMM and the fast direct solver. It is also to be noted that whenever the iterative solver converges, the iterative solver is faster than the fast direct solver. This can be observed from Tables~\ref{table:gaussianTable},~\ref{table:multipleGaussiansTable},~\ref{table:cavityTable} and Figure~\ref{IterationTime40}. It is to be observed that for Gaussian contrast with $\kappa=40$ and $N=305568$, the iterative solver is more than $100$ times faster than the direct solver.

If one were to compute for multiple right-hand sides, we have the following inferences made from the values of T\textsubscript{Hs}, T\textsubscript{GMRES}, and T\textsubscript{Ps} of Tables~\ref{table:gaussianTable},~\ref{table:multipleGaussiansTable},~\ref{table:cavityTable}, and~\ref{table:lensTable}.
\vspace{-0.5pc}
\begin{enumerate}
  \itemsep0em
  \item
  Low wavenumbers: For small $N$, the HODLR direct solver is advantageous over iterative solvers. For large $N$, the iterative solvers are advantageous over the HODLR direct solver.
  \item
   High wavenumbers ($\kappa = 300$): The HODLR direct solver is advantageous over iterative solvers.
\end{enumerate}
The above claims are made with respect to the experiments considered in this article.

On the time complexity front, the claims we made in Table~\ref{table:complexities} are substantiated by Figure~\ref{fig:Scaling}.

\section{Conclusions}
\label{section:Conclusions}
We propose an algebraic fast multilevel summation (DAFMM). Using this we developed an iterative solver for the Lippmann-Schwinger equation in 2D, where all matrix-vector products were computed using DAFMM. The attractive features of DAFMM are:
(i) Low-rank bases are obtained using our new Nested Cross Approximation, (ii) the pivots are picked efficiently in a nested fashion, (iii) Low-rank compressions are problem and domain-specific.
We also present a comparative study of a fast direct solver and the proposed fast iterative solver for the Lippmann-Schwinger equation. In the process, we also propose an efficient preconditioner based on HODLR to further accelerate the fast iterative solver and solve the problems that the iterative solver can not.
In this article, DAFMM is presented for scattering in 2D. Nonetheless, it can be adapted to scattering in 3D as well since i) the admissibility condition for the 3D Helmholtz kernel is also directional~\cite{engquist2010fast} ii) our low-rank construction is algebraic.
In the spirit of reproducible computational science, the implementation of the algorithms developed in this article is made available at \url{https://github.com/vaishna77/Lippmann_Schwinger_Solver}.

\section*{Acknowledgments}
Vaishnavi Gujjula acknowledges the support of Women Leading IITM (India) 2022 in Mathematics (SB22230053MAIITM008880). Sivaram Ambikasaran acknowledges the support of Young Scientist Research Award from Board of Research in Nuclear Sciences, Department of Atomic Energy, India (No.34/20/03/2017-BRNS/34278) and MATRICS grant from Science and Engineering Research Board, India (Sanction number: \\MTR/2019/001241).

\clearpage

\end{document}